\documentclass{elsarticle}

\usepackage{tikz}
\usetikzlibrary{patterns}
\usetikzlibrary{calc}
\usetikzlibrary{shapes.misc}
\tikzset{cross/.style={cross out, draw, 
         minimum size=2*(#1-\pgflinewidth), 
         inner sep=0pt, outer sep=0pt}}

\usepackage{pgfplots}
\pgfplotsset{compat=1.14}
\usepgfplotslibrary{fillbetween}
\pgfdeclarelayer{ft}
\pgfdeclarelayer{bg}
\pgfsetlayers{bg,main,ft}

\usepackage{caption}
\usepackage{subcaption}
\usepackage{bm}
\usepackage{amssymb}
\usepackage{amsmath}
\usepackage{graphicx}
\usepackage{booktabs}
\usepackage{comment}
\usepackage[autolanguage]{numprint}
\usepackage{siunitx}
\usepackage{layouts}
\usepackage{float}
\usepackage{esdiff}
\usepackage{hyperref}
\usepackage{breakurl}

\journal{Comput. Meth. Appl. Mech. Eng.}

\begin{document}

\begin{frontmatter}

\title{A comparison of interpolation techniques for non-conformal high-order discontinuous Galerkin methods}

\author[cemps]{Edward Laughton\corref{mycorrespondingauthor}}
\cortext[mycorrespondingauthor]{Corresponding author}
\ead{el326@exeter.ac.uk}

\author[cemps]{Gavin Tabor}
\ead{g.r.tabor@exeter.ac.uk}

\author[cemps]{David Moxey}
\ead{d.moxey@exeter.ac.uk}

\address[cemps]{College of Engineering, Mathematics and Physical Sciences, University of Exeter, UK.}

\begin{abstract}
The capability to incorporate moving geometric features within models for complex simulations is a common requirement in many fields.  Fluid mechanics within aeronautical applications, for example, routinely feature rotating (e.g. turbines, wheels and fan blades) or sliding components (e.g. in compressor or turbine cascade simulations). With an increasing trend towards the high-fidelity modelling of these cases, in particular combined with the use of high-order discontinuous Galerkin methods, there is therefore a requirement to understand how different numerical treatments of the interfaces between the static mesh and the sliding/rotating part impact on overall solution quality. In this article, we compare two different approaches to handle this non-conformal interface. The first is the so-called mortar approach, where flux integrals along edges are split according to the positioning of the non-conformal grid. The second is a less-documented point-to-point interpolation method, where the interior and exterior quantities for flux evaluations are interpolated from elements lying on the opposing side of the interface. Although the mortar approach has significant advantages in terms of its numerical properties, in that it preserves the local conservation properties of DG methods, in the context of complex 3D meshes it poses notable implementation difficulties which the point-to-point method handles more readily. In this paper we examine the numerical properties of each method, focusing not only on observing convergence orders for smooth solutions, but also how each method performs in under-resolved simulations of linear and nonlinear hyperbolic problems, to inform the use of these methods in implicit large-eddy simulations.
\end{abstract}


\end{frontmatter}

\section{Introduction}

Problems containing features that move or deform are found in many research areas, but are particularly prevalent in the study of various fluid dynamics phenomena. In particular, aeronautical applications commonly feature rotating or sliding geometries, with typical examples in this area including turbomachinery~\cite{tyacke2019,johnston1997}, insect and avian flight aerodynamics~\cite{sun2014, sane2003}, unmanned aerial vehicles~\cite{mi2019,qiwei2020}, and HVAC (heating, ventilating, and air conditioning)~\cite{adam2010,lien2010}. Being able to accurately model these moving geometries and their subsequent impact on the underlying flow physics is highly important: for example, predicting how the profiles of turbine or compressor blades impact on propulsion efficiency, or how wing profiles affect the performance of wind turbines. Moreover, the cost and difficulty of performing full-scale experimental testing of such geometries can be challenging from the perspective of both instrumentation and expense. For these reasons, computational fluid dynamics (CFD) is now commonplace in the design and modelling process. If a CFD method is to be regarded as universally useful in these application areas, then it must be capable of accurately modelling moving geometry and ideally provide high-fidelity results beyond the scope of physical field tests alone.
 
Most leading software for CFD is based around lower-order finite volume or finite element methods, typically leveraging the computationally-cheap Reynolds Averaged Navier-Stokes (RANS) equations in combination with a turbulence closure model. However this approach has natural limitations in studying the aforementioned problems at very high levels of fidelity~\cite{spalart2010}. With the large increases in computational power in recent years, a more recent trend is to instead consider transient simulations that leverage implicit large-eddy simulation (iLES) or under-resolved direct numerical simulation (uDNS)~\cite{aspden2008, bosshard2015}. This approach is more computationally expensive than RANS, but also provides greater accuracy and enables high-fidelity simulations of the complex geometries that lie in this regime~\cite{Chew2018}. The combination of LES with less common high-order methods, either based on continuous or discontinuous Galerkin (CG/DG) methods, has seen significant interest in recent years, particularly in aeronautics applications~\cite{lombard2016}. From a numerical perspective, high-order methods possess far lower levels of numerical diffusion and dispersion, making them ideally suited to resolving features across long time- and length-scales. This can overcome a significant bottleneck when considering these simulations at lower orders, since very fine grid resolutions are required to overcome the effects of numerical errors~\cite{Ghosal1996}. Additionally, from a computational perspective, the larger number of floating-point operations that are required per degree-of-freedom as the polynomial order is increased means that, when equipped with tensor-contraction techniques such as sum factorisation, high-order methods can be used to overcome the memory bandwidth bottlenecks that are common in modern computational hardware~\cite{moxey-2020b}. The combination of these effects means that high-order methods can achieve higher accuracy per degree-of-freedom at equivalent or lower computational cost to lower-order methods.

However, these methods are somewhat less well-explored in the simulation of problems involving rotating or sliding geometries, which require the treatment of non-conformal interfaces between elements. In this article, we explore two common approaches to the handling of non-conformal interfaces and compare their numerical performance in a range of linear and nonlinear problems.

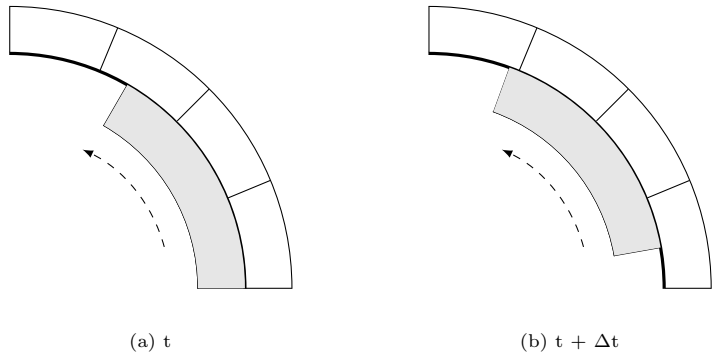
\begin{figure}
\begin{center}
\begin{subfigure}[b]{0.45\textwidth}
\begin{center}
\begin{tikzpicture}[scale=2.5]

\coordinate (A) at (-3,0);
\draw [rotate around={15:(A)}, dashed, -latex] (-2.15,0) arc(0:50:0.8) ;
\begin{scope}

\draw [name path = A, line cap = rect] (-2.0,0) arc(0:60:1);
\draw [name path = AA, line cap = rect] (-1.75,0) arc(0:60:1.25);
\tikzfillbetween[of=A and AA, on layer=ft]{black!10};

\draw [line cap = rect] (-2.0,0) arc(0:60:1)
node[minimum size=0cm, pos=0,inner sep=0pt] (a0) {}
node[minimum size=0cm,pos=0.25,inner sep=0pt] (a1) {}
node[minimum size=0cm,pos=0.5, inner sep=0pt] (a2) {}
node[minimum size=0cm,pos=0.75,inner sep=0pt] (a3) {}
node[minimum size=0cm,pos=1,inner sep=0pt] (a4) {};

\draw [line cap = rect] (-1.75,0) arc(0:60:1.25)
node[minimum size=0cm,pos=0,inner sep=0pt] (aa0) {}
node[minimum size=0cm,pos=0.25,inner sep=0pt] (aa1) {}
node[minimum size=0cm,pos=0.5, inner sep=0pt] (aa2) {}
node[minimum size=0cm,pos=0.75,inner sep=0pt] (aa3) {}
node[minimum size=0cm,pos=1,inner sep=0pt] (aa4) {};

\foreach \t in {0,1,...,4}
\draw [line cap=rect] (a\t) -- (aa\t);
\end{scope}

\begin{scope}

\draw  [very thick, shorten <= -0.2pt, shorten >= -0.25pt] (-1.75,0) arc(0:90:1.25)
node[minimum size=0cm,pos=0,inner sep=0pt] (c0) {}
node[minimum size=0cm,pos=0.25,inner sep=0pt] (c1) {}
node[minimum size=0cm,pos=0.5, inner sep=0pt] (c2) {}
node[minimum size=0cm,pos=0.75,inner sep=0pt] (c3) {}
node[minimum size=0cm,pos=1,inner sep=0pt] (c4) {};

\draw [line cap = rect](-1.5,0) arc(0:90:1.5)
node[minimum size=0cm,pos=0,inner sep=0pt] (b0) {}
node[minimum size=0cm,pos=0.25,inner sep=0pt] (b1) {}
node[minimum size=0cm,pos=0.5, inner sep=0pt] (b2) {}
node[minimum size=0cm,pos=0.75,inner sep=0pt] (b3) {}
node[minimum size=0cm,pos=1,inner sep=0pt] (b4) {};

\foreach \t in {0,1,...,4}
\draw [line cap=rect] (b\t) -- (c\t);
\end{scope}

\end{tikzpicture}
\end{center}
\caption{t}
\end{subfigure}
\begin{subfigure}[b]{0.45\textwidth}
\begin{center}
\begin{tikzpicture}[scale=2.5]

\coordinate (A) at (-3,0);
\draw [rotate around={15:(A)}, dashed, -latex] (-2.15,0) arc(0:50:0.8);

\begin{scope} [rotate around={10:(A)}]

\draw [name path = A, line cap = rect] (-2.0,0) arc(0:60:1);
\draw [name path = AA, line cap = rect] (-1.75,0) arc(0:60:1.25);
\tikzfillbetween[of=A and AA, on layer=ft]{black!10};

\draw [name path = A, line cap = rect] (-2.0,0) arc(0:60:1)
node[minimum size=0cm,pos=0,inner sep=0pt] (a0) {}
node[minimum size=0cm,pos=0.25,inner sep=0pt] (a1) {}
node[minimum size=0cm,pos=0.5, inner sep=0pt] (a2) {}
node[minimum size=0cm,pos=0.75,inner sep=0pt] (a3) {}
node[minimum size=0cm,pos=1,inner sep=0pt] (a4) {};

\draw [name path = AA, line cap = rect] (-1.75,0) arc(0:60:1.25)
node[minimum size=0cm,pos=0,inner sep=0pt] (aa0) {}
node[minimum size=0cm,pos=0.25,inner sep=0pt] (aa1) {}
node[minimum size=0cm,pos=0.5, inner sep=0pt] (aa2) {}
node[minimum size=0cm,pos=0.75,inner sep=0pt] (aa3) {}
node[minimum size=0cm,pos=1,inner sep=0pt] (aa4) {};

\foreach \t in {0,1,...,4}
\draw [line cap=rect] (a\t) -- (aa\t);
\end{scope}

\begin{scope} 
\draw [very thick, shorten <= -0.2pt, shorten >= -0.25pt] (-1.75,0) arc(0:90:1.25)
node[minimum size=0cm,pos=0,inner sep=0pt] (c0) {}
node[minimum size=0cm,pos=0.25,inner sep=0pt] (c1) {}
node[minimum size=0cm,pos=0.5, inner sep=0pt] (c2) {}
node[minimum size=0cm,pos=0.75,inner sep=0pt] (c3) {}
node[minimum size=0cm,pos=1,inner sep=0pt] (c4) {};

\draw [line cap = rect] (-1.5,0) arc(0:90:1.5)
node[minimum size=0cm,pos=0,inner sep=0pt] (b0) {}
node[minimum size=0cm,pos=0.25,inner sep=0pt] (b1) {}
node[minimum size=0cm,pos=0.5, inner sep=0pt] (b2) {}
node[minimum size=0cm,pos=0.75,inner sep=0pt] (b3) {}
node[minimum size=0cm,pos=1,inner sep=0pt] (b4) {};

\foreach \t in {0,1,...,4}
\draw [line cap=rect] (b\t) -- (c\t);
\end{scope}

\end{tikzpicture}
\end{center}
\caption{t + $\Delta$t}
\end{subfigure}

\end{center}
\caption{An example domain showing an inner rotating region and an outer stationary region bounded by a non-conforming interface zone, showing the differing physical positions at two different times.}
\label{fig:nonconformal}
\end{figure}

\subsection{Requirements for moving geometry simulations}

One approach to tackling the problem of moving geometry is the sliding mesh method, where the mesh is separated into two or more separate regions, and during the simulation the regions will slide relative to one another. This provides a way to prescribe simple mesh motion via rotation or translation. The most basic problem case is to employ a stationary outer region, with a rotating circular region within it, which is found in many applications, for example modelling flow in a stirred tank~\cite{Bakker1997}. An exaggerated example of this arrangement is shown in figure~\ref{fig:nonconformal}, where it is clear that this process results in a non-conformal mesh: i.e. a mesh where elements do not connect to precisely one other element through one of their edges or faces. Most CFD simulations make use of {\em conformal} meshes, where each edge (in 2D) or face (in 3D) of an element has precisely one neighbouring element. As such, techniques need to be developed in order to accurately preserve solution quality across the non-conformal interface.

\begin{figure}
\begin{center}
\begin{tikzpicture}[scale=5]

  \draw[fill=black!10] (-0.4,0.66) -- (0.3, 0.66) -- (0.3, 1) -- (-0.4,1);
  \draw[fill=black!10] (-0.4,0.33) -- (0.3, 0.33) -- (0.3, 0.66) -- (-0.4,0.66);
  \draw[fill=black!10] (-0.4,0) -- (0.3, 0) -- (0.3, 0.33) -- (-0.4,0.33);

  \draw[fill=black!10]  (1.9, 1) --(1.2,1) -- (1.2,0.5) -- (1.9, 0.5) ;
	  \draw[fill=black!10]  (1.9, 0.5) --(1.2,0.5) -- (1.2,0) -- (1.9, 0) ;

\draw (0.75, 0) -- (0.75, 1);
	   
  \draw (-0.015, 0.825) node{Element A, $\Omega_{A}$};
  \draw (-0.015, 0.495) node{Element B, $\Omega_{B}$};
    \draw (-0.015, 0.165) node{Element C, $\Omega_{C}$};

  \draw (1.55, 0.75) node{Element D,  $\Omega_{D}$};
  \draw (1.55, 0.25) node{Element E,  $\Omega_{E}$};

  \draw[fill=black] (0.75, 0) circle (0.5pt);
  \draw[fill=black] (0.75, 0.33) circle (0.5pt);
    \draw[fill=black] (0.75, 0.50) circle (0.5pt);
  \draw[fill=black] (0.75, 0.66) circle (0.5pt);
  \draw[fill=black] (0.75, 1) circle (0.5pt);
   \draw[fill=black] (0.75, 1)  circle (0.5pt) node[above,yshift=3pt]{Mortars};

      \draw[dotted,thick,latex-latex, shorten >= 3pt, shorten <= 3pt]  (0.3, 0.825) -- (0.75,0.825);
      \draw[dotted,thick,latex-latex, shorten >= 3pt, shorten <= 3pt]  (0.3, 0.5) -- (0.75,0.58);
      \draw[dotted,thick,latex-latex, shorten >= 3pt, shorten <= 3pt]  (0.3, 0.5) -- (0.75,0.42);
      \draw[dotted,thick,latex-latex, shorten >= 3pt, shorten <= 3pt]  (0.3, 0.165) -- (0.75,0.165);
      
      \draw[dotted,thick, latex-latex, shorten <= 3pt, shorten >= 3pt]  (1.2,0.75) -- (0.75, 0.825) ;
      \draw[dotted,thick, latex-latex, shorten <= 3pt, shorten >= 3pt]  (1.2,0.75) -- (0.75,0.58) ;
      \draw[dotted,thick, latex-latex, shorten <= 3pt, shorten >= 3pt]  (1.2,0.25) -- (0.75,0.42);
      \draw[dotted,thick, latex-latex, shorten <= 3pt, shorten >= 3pt]  (1.2,0.25) -- (0.75, 0.165); 

\end{tikzpicture}
\end{center}
\caption{Mortar construction for a non-conformal interface showing the connection between elements and mortars.}
\label{fig:mortars}
\end{figure}
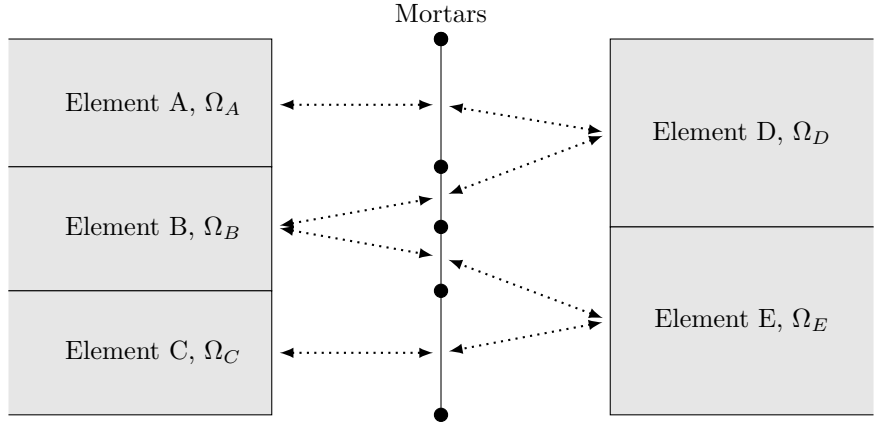

In the `classical' spectral element method, where $C^0$ continuity is imposed between elements in a CG formulation, three main techniques have been evaluated for use in non-conformal meshes. Possibly the most well-known of these arises when performing $h$-adaptation in an octree-like manner, so that 2-to-1 element subdivisions are obtained in the resulting mesh. In this case, hanging nodes are generated and their values can be constrained through the analytic definition of the basis functions lying along an edge or face, together with the assembly mapping that is used to construct mass and stiffness matrices~\cite{demkowicz1989, offermans2019, bangerth2007}. However, in the sliding mesh case, elements may overlap at arbitrary positions along their edges and faces, making this approach infeasible. Possibly the most widely-adopted approach to implementing generic non-conformal interfaces in the CG setting is the mortar technique~\cite{mavriplis}. In this setting, one augments the traditional $C^0$ function spaces for each conformal domain with functions defined on mortar elements at the interface between two domains. The weak form of the problem is then augmented to incorporate a penalty for the jump across the interface in an appropriate manner, so that the convergence order of the scheme is retained. This approach is visualised in figure~\ref{fig:mortars}, where we note that the mortar elements are constructed at the common intersection points of each element.

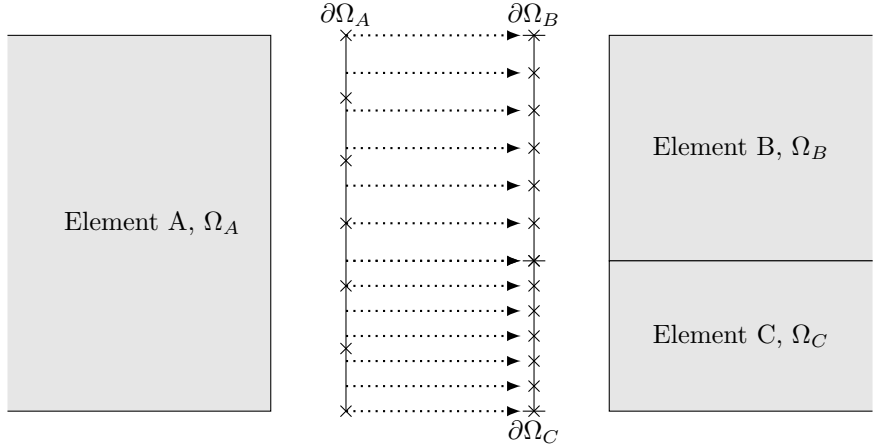
\begin{figure}
\begin{center}
\begin{tikzpicture}[scale=5]

    \draw[fill=black!10] (-0.4,0) -- (0.3, 0) -- (0.3, 1) -- (-0.4,1);
    \draw[fill=black!10]  (1.9, 1) --(1.2,1) -- (1.2,0.4) -- (1.9, 0.4) ;
    \draw[fill=black!10]  (1.9, 0) --(1.2,0) -- (1.2,0.4) -- (1.9, 0.4) ;
  
    \draw (0.5, 0) -- (0.5, 1); 
    \draw (1,1) -- (1,0);
    
    \draw (0.97,0.0) -- (1.03,0.0);
    \draw (0.97,0.4) -- (1.03,0.4);
    \draw (0.97,1.0) -- (1.03,1.0);
  
    \draw (-0.015, 0.5) node{Element A, $\Omega_{A}$};
    \draw (1.55, 0.7) node{Element B,  $\Omega_{B}$};
    \draw (1.55, 0.2) node{Element C,  $\Omega_{C}$};

    \foreach \x in {0,...,6} {
        \draw[dotted,thick,-latex, shorten >= 5pt]  (0.5, 0.4+\x/10) -- (1,0.4+\x/10);
        \draw[dotted,thick,-latex, shorten >= 5pt]  (0.5, 0.4*\x/6) -- (1,0.4*\x/6);

        \draw (1, 0.4+\x/10) node[cross=2.5pt] {};
        \draw (0.5, \x/6) node[cross=2.5pt] {};
        \draw (1, 0.4*\x/6) node[cross=2.5pt] {};
    }
    
    \node[above] at (0.5, 1) {$\partial\Omega_{A}$};
    \node[above] at (1, 1.0) {$\partial\Omega_{B}$};
    \node[below] at (1, 0.0) {$\partial\Omega_{C}$};
\end{tikzpicture}
\end{center}
\caption{Schematic representation of point-to-point interpolation across a non-conformal interface. Crosses represent the integration points on the respective trace of each element. Dotted arrows denote a high-order evaluation of the polynomial on the edge of element A at the points required for flux evaluations of elements B and C.}
\label{fig:interp}
\end{figure}

An alternative approach to imposing non-conformal conditions, and one which is perhaps less commonly-used, is to instead adopt a point-to-point interpolation across the elemental interface. In this setting, no attempt to construct mortar elements is made and the function space is defined in the usual manner for each conformal domain. However, when values within elements are desired at the left-hand side of the interface, they are obtained by performing a polynomial interpolation from the values on the right-hand side, and vice versa. This approach was first implemented and tested for geophysical problems in~\cite{rosenberg-2006} in the CG setting, where it was shown to demonstrate convergence-order preserving properties. Within the DG setting this allows for the the values of the exterior conserved variables to be obtained from interpolating the interior trace values on the opposing side of the interface. A sample visualisation of this approach is shown in figure~\ref{fig:interp}, where dotted arrows denote the evaluation of the high-order polynomial defined by points on the edge of element $\Omega_A$ to obtain their values within the boundaries of elements in $\partial\Omega_B$ and $\partial\Omega_C$. This interpolation process can be built into the assembly operation that is used to construct mass and stiffness matrices.

\subsection{Non-conformal techniques for the discontinuous Galerkin method}

At present, there is a significant interest in the development of high-order fluid dynamics solvers for iLES/uDNS based around the DG method due to its favourable stability properties in these settings~\cite{moura2019,moura2017,flad2017}. In the context of sliding mesh simulations, DG also offers an easier route to the accurate treatment of non-conformal interfaces between elements across the sliding interface, since elements are naturally disconnected as part of the formulation of the method. Additionally, approaches to impose non-conformal interfaces in DG are perhaps less well-explored than in CG.

In the DG formulation, connectivity between elements is imposed through a flux term, which may be either an upwind-type solver for linear problems or a more complex Riemann problem for more general nonlinear hyperbolic systems. These fluxes are computed on integrals across each edge of an element and take the form
\[
\int_{\Gamma_e} \tilde{\bm{f}}(\bm{u}^+, \bm{u}^-)\cdot\bm{n}\, ds
\]
where $\Gamma_e$ is an edge of element $\Omega^e$, $\bm{u}^+$ and $\bm{u}^-$ is a vector of conserved variables on the exterior and interior of the element respectively, $\tilde{\bm{f}}$ is the numerical flux function and $\bm{n}$ is an outwards-facing normal. The question then is how one computes these integrals, given that the exterior values $\bm{u}^+$ may now lie across more than one element on the other side of the interface.

The mortaring approach has been investigated in a number of works from the DG perspective. First, we note that unlike the CG setting which requires modifications to the function space and weak form of the problem, in DG by `mortaring' we only refer to the act of constructing mortar elements on which to compute the flux integral. That is, the integral above is split into multiple integrals, one for each mortar element. This approach was first investigated by Kopriva et al. in the study of both fluid dynamics~\cite{kopriva} and electromagnetics problems~\cite{kopriva2}. The same approach has been used fairly extensively for problems involving sliding meshes; for example by~\cite{ferrer-2012} in a hybrid DG-Fourier pseudospectral solver for the incompressible Navier-Stokes equations, in the construction of a spectral difference solver for the compressible Navier-Stokes equations that incorporates sliding grids by Zhang and Liang~\cite{Zhang2015} and more recently in the hyperbolic solver FLEXI~\cite{krais-2020}. This approach has the significant advantage that it preserves the local conservation property of the DG method, which is important from the perspective of obtaining accurate results that conserve mass (in the case of CFD). However, although mortaring is straightforward in two dimensions, a significant challenge in the use of the mortar approach for general three-dimensional problems is the generation of the mortar elements themselves. In general this could involve the re-meshing of the non-conformal interface between domains at each timestep in order to generate an appropriate mortar space, as adopted by Aguerre et al.~\cite{aguerre-2017} for finite volume simulations based on the supermesh construction of Farrell et al.~\cite{farrell-2009}. We note that these works consider only straight-sided elements. In some settings, such as a sliding mesh defined by a translation, this approach could therefore readily be adopted to a high-order setting. However, in problems involving rotation such as in figure~\ref{fig:nonconformal}, the interface between regions now becomes curved, thereby significantly increasing the complexity and computational requirements in this approach.

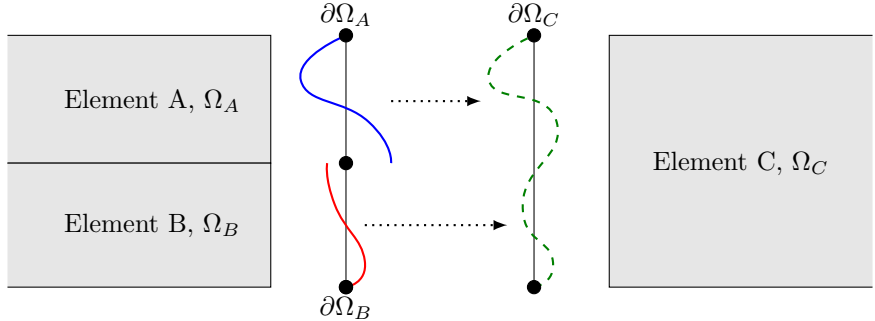
\begin{figure}
\begin{center}
\begin{tikzpicture}[scale=5]

  \draw[fill=black!10] (-0.4,0.66) -- (0.3, 0.66) -- (0.3, 1) -- (-0.4,1);
  \draw[fill=black!10] (-0.4,0.33) -- (0.3, 0.33) -- (0.3, 0.66) -- (-0.4,0.66);

  \draw[fill=black!10]  (1.9, 1) --(1.2,1) -- (1.2,0.33) -- (1.9, 0.33) ;
	
\draw (0.5, 0.33) -- (0.5, 1);
 \draw[blue, thick] plot [smooth, tension=1] coordinates {(0.62, 0.66) (0.56, 0.77) (0.38,0.88) (0.5, 1)} ;
  \draw[red, thick] plot [smooth, tension=1] coordinates {(0.5, 0.33) (0.55, 0.40) (0.47,0.55) (0.45, 0.66)} ;

\draw (1,1) -- (1,0.33);	
    \draw[green!50!black, thick, dashed] plot [smooth, tension=1] coordinates {(1, 0.33) (1.05, 0.40) (0.97,0.55) (1.06, 0.77) (0.88,0.88) (1, 1)} ;
	   
  \draw (-0.015, 0.825) node{Element A, $\Omega_{A}$};
  \draw (-0.015, 0.495) node{Element B, $\Omega_{B}$};
  \draw (1.55, 0.66) node{Element C,  $\Omega_{C}$};

  \draw[fill=black] (0.5, 0.33) circle (0.5pt);
   \draw[fill=black] (0.5, 0.66) circle (0.5pt);
   \draw[fill=black] (0.5, 1)  circle (0.5pt) node[above]{$\partial\Omega_{A}$};
      \draw[fill=black] (0.5, 0.33)  circle (0.5pt) node[below]{$\partial\Omega_{B}$};
   
   \draw[fill=black] (1, 0.33) circle (0.5pt);
   \draw[fill=black] (1, 1.0)  circle (0.5pt) node[above]{$\partial\Omega_{C}$};
   
      \draw[dotted,thick,-latex, shorten >= 3pt]  (0.62, 0.825) -- (0.88,0.825);
      \draw[dotted,thick,-latex, shorten >= 3pt]  (0.55, 0.495) -- (0.95,0.495);
  
\end{tikzpicture}
\end{center}
\caption{Interpolation across a non-conformal 2D interface showing the solution of element A (\textbf{\color{blue} ---}), element B (\textbf{\color{red} ---}) and the resultant discontinuity in the solution of Element C (\textbf{\color{green!50!black} -- --}), with element boundaries indicated by the circles ($\bullet$).}
\label{fig:cont}
\end{figure}

The alternative approach is therefore to consider the point-to-point interpolation method in the DG context, since implementation is relatively straightforward by comparison as it does not require the construction of mortar elements. However, neither the implementation, performance or  robustness of this approach for has been thoroughly investigated in the literature to date. In particular, potential issues may arise from the discontinuity of fluxes between elements: for smooth solutions and at high polynomial orders, the interpolation between neighbouring non-conformal elements will likely introduce very little error into the resulting solution. However, in the presence of under-resolved simulations, which are more prone to admitting discontinuities in flow solution between elements, the discontinuity may introduce additional numerical error that warrants further study. A prototypical example which demonstrates this in an illustrative manner is shown in Fig,~\ref{fig:cont}. On the left side of the interface, the two discontinuous solutions from elements $\Omega_{A}$ and $\Omega_{B}$ must be sampled at integration points on the skeleton of $\Omega_{C}$. If the two functions are sufficiently discontinuous, the interpolation procedure could result in spurious noise introduced into the interior of $\Omega_{C}$.

\subsection{Aim of this work}

To date, the point-to-point method has not been well-studied in the literature.  A study by Kopera and Giraldo~\cite{kopera-2015} is one of the very few references, to the best of the authors' knowledge, that consider the point-to-point interpolation approach in DG, where CG and DG implementations of the interpolation technique are examined and their mass conservation properties are reported. However we note that in this case, only hanging-node type vs. more generic non-conformal interfaces are considered. Additionally, this work was performed in well-resolved cases which may not be the case for more general iLES/uDNS-type problems. In this article, we therefore aim to address this gap in the understanding of the performance of these approaches by performing a comparative study of the mortar and point-to-point techniques. We consider several aspects, including a validation of convergence order for both approaches, the performance of each method in terms of numerical diffusion for a linear transport equation at varying degrees of underresolution, the behaviour of each method when considering the nonlinear problem of the compressible Euler equations across long time periods, and the compressible Navier-Stokes equations in a 3D setting.

The remainder of the paper is structured as follows. In section~\ref{sec:theory} we set out the theoretical framework of the two formulations and outline our implementation strategy within the spectral/$hp$ element framework Nektar++~\cite{Cantwell2015,moxey-2020a}. Section~\ref{sec:results} presents the results of our studies for a linear transport equation and the nonlinear compressible Euler equations in two dimensions. In section~\ref{sec:3d}, we consider more realistic fluid simulations in three dimensions, by examining the transition to turbulence in a canonical Taylor-Green vortex problem. Finally, in section~\ref{sec:conclusions}, we draw some brief conclusions and discuss the key performance characteristics of each method.

\section{Theory}
\label{sec:theory}

\subsection{The DG formulation of the spectral/\textit{hp} element method}

In this section, we briefly introduce the discontinuous Galerkin (DG) discretisation of the spectral/\textit{hp} element method. A more thorough overview can be found in several other works, e.g.~\cite{Cantwell2015,moxey-2020a,karniadakis2005}. The starting point for the DG formulation is the same as any other typical finite element problem, in that we consider a domain $\Omega$ comprised of non-overlapping elements $\Omega^e$ such that $\Omega = \bigcup_e \Omega^e$. Given a general hyperbolic conservation law for conserved variables $\bm{u}$ taking the form
\begin{equation}
\diffp{{\bm{u}}}{t} + \nabla\cdot\bm{F}(\bm{u}) = 0,
\label{eq:cons}
\end{equation}
we follow the standard Galerkin approach and, on a single element, construct the weak form via multiplication by a test function $v$ and integrating by parts to obtain
\begin{equation}
  \left(v, \diffp{{\bm{u}}}{t}\right)_{\Omega^e} + \left\langle v\bm{n}, \tilde{\bm{f}}(\bm{u}^+, \bm{u}^-) \right\rangle_{\partial\Omega^e}
  - \left(\nabla v, \bm{F}(\bm{u}) \right)_{\Omega^e} = 0,
  \label{eq:weak}
\end{equation}
where $(u,v)_{\Omega^e} = \int_{\Omega^e} uv\, d\bm{x}$ and $\langle u, v\rangle_{\partial\Omega^e}= \int_{\partial\Omega^e} uv\, ds$ denote inner products on the volume and surface, respectively. Moreover, $\tilde{\bm{f}}$ defines a numerically-calculated flux term which, as explained in the previous section, may take the form of a general Riemann problem, and which depends on the element-exterior and interior velocities $\bm{u}^+$ and $\bm{u}^-$, respectively. Within each element, we represent $\bm{u}$ using an expansion of high-order polynomials, so that
\[
\bm{u}^\delta = \sum_n \hat{\bm{u}}_n\phi_n\left([\bm{\chi}^e]^{-1}(\bm{x})\right).
\]
In this expression, we note that the approximation is defined with the use of a standard (reference) element $\Omega_{\text{st}}$, with $\phi_n$ denoting an appropriate set of basis functions. An isoparametric mapping $\bm{\chi}^e:\Omega_{\text{st}}\to\Omega^e$ defines a possibly curvilinear element $\Omega^e$, so that $\bm{x}=\bm{\chi}^e(\bm{\xi})$ for $\bm{\xi}\in\Omega_{\text{st}}$. We additionally equip the standard element with a distribution of quadrature points $\bm{\xi}_q$ and weight $w_q$, so that upon selecting test functions $v = \phi_n$ we then evaluate the terms in eq.~\eqref{eq:weak} as finite summations, i.e.
\[
\left(\nabla\phi_n, \bm{F}(\bm{u})\right)_{\Omega^e} \approx \sum_{q} \nabla\bm{\chi}^e(\bm{\xi}_q)^{-T}\nabla\phi_n(\bm{\xi}_q)\cdot\bm{F}(\bm{u}(\bm{x}_q)) \det\left(\bm{\chi}^e(\bm{\xi}_q)\right) w_q
\]
In this study, we consider only two-dimensional elements and select tensor products of Gauss-Lobatto-Legendre points to evaluate quadrature. As basis functions we adopt the hierarchical modified basis of Karniadakis \& Sherwin~\cite{karniadakis2005}. Similarly to the classical Lagrange basis, these basis functions have the beneficial property of boundary-interior decomposition, which makes the addition of flux terms into the overall elemental degrees of freedom a straightforward addition operation. In particular we note that the flux integral terms can be considered along each edge $i$ of $\Omega^e$, which we denote by $\Gamma^e_i$, as the integral
\[
\int_{\Gamma^e_i} \psi_n(\bm{\xi}) \tilde{\bm{f}}(\bm{u}^+, \bm{u}^-) ds
\]
where now $\psi_n$ denotes a basis function with support along edge $i$. In particular, we note that the solution variables along $\bm{y}\in\Gamma^e_i$ can be written as a polynomial expansion
\begin{equation}
    \bm{u}(\bm{y}) = \sum_n \hat{\bm{u}}_n \psi_n(\bm{y})
    \label{eq:polyinterp}
\end{equation}
As alluded to in the introduction, the central focus of this work is to understand how different evaluations of the flux term in the presence of a non-conformal mesh influence the overall properties and stability of the DG method. In the following sections, we outline the formulation of both the point-to-point interpolation method and the mortar method.

\subsection{The point-to-point interpolation method}

In the point-to-point interpolation method, the interface is handled using a direct interpolation from one side to the other. That is, when we require the values of exterior conserved variables $\bm{u}^+$ at a spatial position $\bm{y}^*\in\Gamma^e_i$, we adopt the following procedure:

\begin{itemize}
    \item determine a corresponding element $\Omega^f$ that contains the point $\bm{y}^*$ along an edge $\Gamma^f_j$;
    \item perform a polynomial interpolation at that position using eq.\eqref{eq:polyinterp} in order to determine $\bm{u}^+$.
\end{itemize}

This interpolation is performed for every integration point along $\Gamma^e_i$, as shown diagrammatically in figure~\ref{fig:interp}. Once the trace space (i.e. the collection of all edges in the interface) has been fully populated by interpolation, the DG solver can continue as usual with a Riemann solver to calculate the numerical flux to then be added into elemental coefficient spaces. 
In order to determine a corresponding element that contains the point $\bm{y}^*$, we require the ability to determine the distance of a desired point from any given edge $\Gamma^f_j$. For edges that are straight-sided, this translates into a simple geometric problem which may be solved analytically. However, for curvilinear elements, we must instead utilise the parametric mapping for the two-dimensional element, which gives a coordinate mapping $\bm{x} = \bm{\chi}^e(\bm{\xi})$. In particular, for each edge in the non-conformal interface, we minimise an objective function $d(\bm{\xi}; \bm{y}^*) = \|\bm{x}-\bm{y}^*\|_2^2 = \|\bm{\chi}^e(\bm{\xi})-\bm{y}^*\|_2^2$, i.e. the square of the Euclidean norm $\| \cdot\|_2$ between a point $\bm{\xi}$ within the edge and the target point $\bm{y}^*$. This then allows us to determine the corresponding reference space point $\bm{\xi}^* = \min_\xi d(\bm{\xi}; \bm{y}^*)$, so that $\bm{\chi}(\bm{\xi})^*$ has minimum distance to $\bm{y}^*$. An edge that has $\bm{\chi}(\bm{\xi}^*)\approx\bm{y}^*$ is chosen as the corresponding edge $\Gamma^f_j$. In our implementation, this is solved via a gradient-descent method utilising a quasi-Newton search direction and backtracking line search, but other Newton-type methods will provide similar convergence properties. Since this is additionally an expensive operation to be performed for every edge within the interface, we make use of an $r$-tree structure to reduce the initial search space. The octants that are used to construct the $r$-tree are defined as the bounding box for each curvilinear edge. In this manner, the $r$-tree can first be interrogated to determine a subset of possible edges under which to then perform the nonlinear optimisation of distance, which further reduces computational cost.

Finally, in the minimisation process above, we require the evaluation of each polynomial expansion~\eqref{eq:polyinterp} at any arbitrary point $\bm{\xi}$ in the reference element. Although this can be computed directly from eq.~\eqref{eq:polyinterp}, this would require the evaluation of each basis function at the same arbitrary point. We note that, for numerical integration purposes within the DG scheme, we already naturally represent $\bm{u}$ at solution points $\bm{y}_q = \bm{\chi}^e(\bm{\xi}_q)$ that correspond with Gauss-Lobatto quadrature points in the reference element $\bm{\xi}_q\in[-1,1]$. This allows us to rewrite eq.~\eqref{eq:polyinterp} as a summation in terms of Lagrange interpolants $\ell_q(\bm{\xi})$ defined using these same points, so that
\[
\bm{u}(\bm{y}) = \sum_q \bm{u}(\bm{y}_q) \ell_q(\bm{\xi}).
\]
Classically, given this collocation representation, one would then generate a diagonal interpolation matrix $\mathcal{I}(\bm{y})$ as outlined in~\cite{karniadakis2005}, and perform a dot product against a vector of points $\bm{u}(\bm{y}_q)$ to evaluate $\bm{u}(\bm{y})$. However, our timings demonstrate that the use of fast summation based on barycentric interpolation techniques described in~\cite{berrut-2004} yield far better performance for this operation. The extension of this technique to higher dimensions is discussed further in section~\ref{sec:3d}.

\subsection{The mortar method}

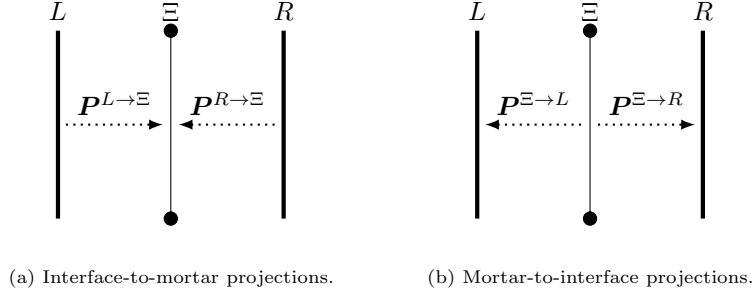
\begin{figure}
\begin{center}
\begin{subfigure}[b]{0.45\textwidth}
\begin{center}
\begin{tikzpicture}[scale=5]

  \draw[fill=black,ultra thick] (0, 0) -- (0, 0.5) node[above]{$L$};
    \draw[fill=black] (0.3, 0) circle (0.5pt) -- (0.3, 0.5) circle (0.5pt) node[above]{$\Xi$};
  \draw[fill=black,ultra thick] (0.6, 0) -- (0.6, 0.5) node[above]{$R$};

         \draw[dotted,thick, -latex, shorten <= 3pt, shorten >= 3pt]  (0,0.25) -- (0.3, 0.25) node[midway, above] {$\bm{P}^{L\rightarrow\Xi}$}; 
         \draw[dotted,thick, latex-, shorten <= 3pt, shorten >= 3pt]  (0.3,0.25) -- (0.6, 0.25) node[midway, above] {$\bm{P}^{R\rightarrow\Xi}$}; 

\end{tikzpicture}
\end{center}
\caption{Interface-to-mortar projections.}
\label{fig:LRmortar}
\end{subfigure}
\begin{subfigure}[b]{0.45\textwidth}
\begin{center}
\begin{tikzpicture}[scale=5]

  \draw[fill=black,ultra thick] (0, 0) -- (0, 0.5) node[above]{$L$};
    \draw[fill=black] (0.3, 0) circle (0.5pt) -- (0.3, 0.5) circle (0.5pt) node[above]{$\Xi$};
  \draw[fill=black,ultra thick] (0.6, 0) -- (0.6, 0.5) node[above]{$R$};

         \draw[dotted,thick, latex-, shorten <= 3pt, shorten >= 3pt]  (0,0.25) -- (0.3, 0.25) node[midway, above] {$\bm{P}^{\Xi\rightarrow L}$}; 
         \draw[dotted,thick, -latex, shorten <= 3pt, shorten >= 3pt]  (0.3,0.25) -- (0.6, 0.25) node[midway, above] {$\bm{P}^{\Xi\rightarrow R}$}; 

\end{tikzpicture}
\end{center}
\caption{Mortar-to-interface projections.}
\label{fig:mortarLR}
\end{subfigure}
\end{center}
\caption{The left, $L$, and right, $R$, interface edges projection relationships from and to the mortar element, $\Xi$.}
\label{fig:LR}
\end{figure}

The second approach we will consider in this paper is the mortar method which maintains the local conservation properties of DG by constructing mortar elements as visualised in figure~\ref{fig:mortars}.  This method applied to the spectral element method was originally developed by Maday et al. \cite{maday}, and has been used for both incompressible flow \cite{mavriplis} and compressible flow problems~\cite{kopriva}. We note again that `mortar' in this sense refers to the act of construction of mortar elements so that flux integrals may be expressed as
\[
\int_{\Gamma^e_i} \psi_n(\bm{\xi}) \tilde{\bm{f}}(\bm{u}^+, \bm{u}^-) ds = \sum_{m=1}^M \int_{\Xi_m} \psi_n(\bm{\xi}) \tilde{\bm{f}}(\bm{u}^+, \bm{u}^-) ds
\]
where $M$ is the number of mortars on edge $i$, and $\Xi_m$ denotes each mortar element. We then construct a polynomial expansion on each mortar element of the same polynomial order. The mortar method is realised by projecting variables from across the interface onto its corresponding mortar element, solving the Riemann problem on the mortars, and then performing an $L^2$ projection in order to consolidate the contributions from each mortar element. The number of mortars connected to a single interface edge and their relative size is arbitrary, allowing for a wide range of varying mesh circumstances. To give a more concrete definition of the method, we utilise the notation prevalent in Zhang and Liang~\cite{Zhang2015} and Kopriva et al.~\cite{kopriva2}, labelling the two contributing interface segments `L' and `R' as shown in figure~\ref{fig:LR}. 

First we recall that each edge in the interface $\Gamma_i^e$ may be represented on a standard segment $-1 \leq \xi \leq 1$ and then mapped using the isoparametric mapping $\bm{\chi}^e$. Similarly, each mortar element has a similar mapping $-1 \leq z \leq 1$ and, in particular, we may write the relationship between the two as
\begin{equation*}
\xi = o + sz
\end{equation*}
where $o$ is the offset of the centre of the mortar relative to the centre of the interface edge, and $s$ is the relative scale factor. The solution on an interface edge can be represented by eq.~\eqref{eq:polyinterp}, so that
\begin{equation*}
u^{\Omega}(\xi) = \sum_p \hat{u}^{\Omega}_{p}\phi_{p}(\xi),
\end{equation*}
where we consider now only a single scalar quantity $u$ for clarity. We can similarly define the solution on the mortar element, $\Xi$, as
\begin{equation*}
u^{\Xi}(z) = \sum_p \hat{u}^{\Xi}_{p}\phi_{p}(z).
\end{equation*}
To project the solutions from the element onto the mortar we minimise the norm in the $L^2$ sense, i.e.
\begin{equation*}
\int^{1}_{-1}{\big(u^{\Xi}(z) -u^{\Omega}(\xi)\big)\phi_{j}(z)\,dz}=0, \quad \text{for all } j.
\end{equation*}
When evaluated at all quadrature points, this can be expressed in matrix form as
\begin{equation*}
\bm{\hat{u}}^{\Xi}= \bm{P}^{\Omega\rightarrow \Xi} \bm{\hat{u}}^{\Omega} =\bm{M}^{-1}\bm{S}^{\Omega\rightarrow\Xi}\bm{\hat{u}}^{\Omega},
\end{equation*}
where $\bm{M}$ is the standard elemental mass matrix, and $\bm{S}^{\Omega\rightarrow\Xi}$ are constructed as
\begin{equation*}
S_{i,j}=\int^{1}_{-1}{\phi_{i}(o + sz)\phi_{j}(z) dz}, \quad \text{for all } i, j.
\end{equation*}

To apply the mortar method to the DG formulation, we therefore adopt the following approach:
\begin{itemize}
    \item Construct both the left and right solutions $\bm{u}^+$ and $\bm{u}^-$ onto the mortar using the projection matrices $\bm{P}^{L\rightarrow \Xi}$ and $\bm{P}^{R\rightarrow \Xi}$ as shown in figure~\ref{fig:LRmortar}.
    \item Once the solutions are on the mortar, the Riemann solver can be used to compute the numerical flux $\tilde{\bm{f}}$.
    \item Projecting from the mortars back onto the interface element requires minimising the trace quantities norm in the L\textsuperscript{2} sense. For $M$ mortars to the interface element $\Omega$ this is as follows
    \begin{equation*}
        \sum_{i = 1}^{M}{\left[ \int_{\Xi_{i}}{\big(\hat{f}^{\Omega}(\xi) -\hat{f}^{\Xi_{i}}(z)\big)\phi_{j}(\xi) d\xi}\right]}  = 0, \quad \text{for all } j.
    \end{equation*}
    In matrix form the solution to this is
    \begin{equation*}
        \bm{\hat{f}}^{\Omega}= \sum_{i=1}^{N}{\left[\bm{P}^{\Xi_{i}\rightarrow \Omega}\bm{\hat{f}}^{\Xi_{i}}\right]} = \sum_{i=1}^{N}{\left[s_{\Xi_{i}}\bm{M}^{-1}\bm{S}^{\Xi_{i}\rightarrow\Omega}\bm{\hat{f}}^{\Xi_{i}}\right]},
    \end{equation*}
    where $\bm{S}^{\Xi\rightarrow\Omega}$ is the transpose of $\bm{S}^{\Omega\rightarrow\Xi}$ taking care to include the respective scale factors.
\end{itemize}
It is also worth noting that for where the geometry of the interface element is identical to the mortar, for example between $\Omega_{A}$ and its corresponding mortar in figure~\ref{fig:mortars}, the projection matrix is merely the identity matrix so $\bm{\hat{u}}^{\Xi} = \bm{\hat{u}}^{\Omega}$ and $\bm{\hat{f}}^{\Omega} = \bm{\hat{f}}^{\Xi}$. This can be used to reduce computational costs.

\section{Results}
\label{sec:results}
In this section, we report on the results of a number of two-dimensional tests using both linear and nonlinear problems to evaluate the efficacy of both the mortar and point-to-point interpolation technique. At each stage, we use conformal grids of similar resolutions to provide a benchmark against which to compare. Each method has been implemented within the Nektar++ spectral/$hp$ element framework~\cite{Cantwell2015, moxey-2020a}. In all cases, we consider only explicit timestepping methods with the use of a standard 4th-order Runge-Kutta time integration scheme unless otherwise stated. The timestep used for each case is reported separately.
\subsection{Convergence order}
\begin{figure}[h]
    \centering
    \includegraphics[width=\textwidth]{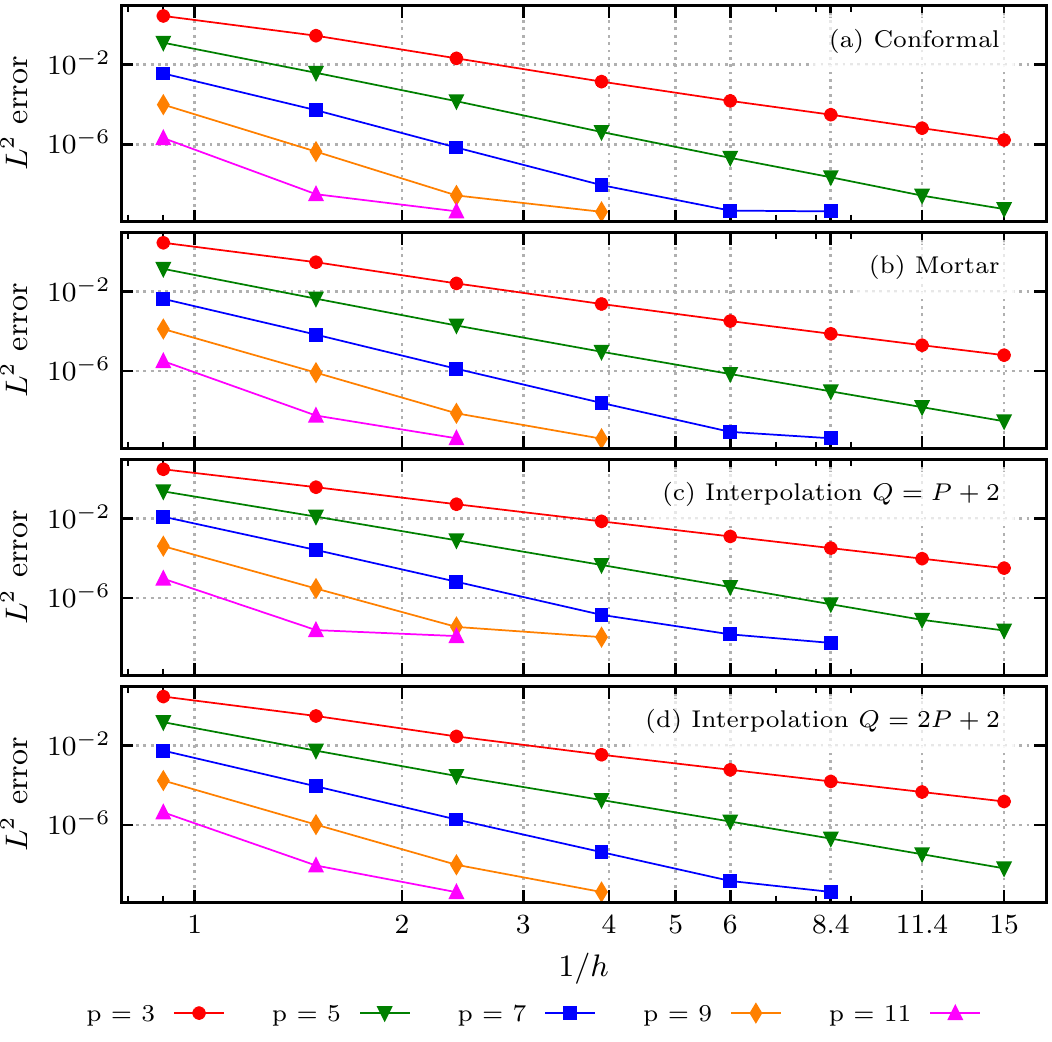}
    \caption{Convergence properties for the conformal case, the non-conformal mortar method, and the non-conformal point-to-point interpolation method.}
    \label{fig:convergence}
\end{figure}
In this first case, we test the correctness of our implementation by performing a standard $h$-convergence study for various polynomial orders $p$. For this, we select a standard linear transport equation within a domain $\Omega=[-5,5]^2$, so that in eq.~\eqref{eq:cons}, $\bm{F}(u) = \bm{v}u$ for a constant velocity $\bm{v} = (1,0)$. We select an initial condition that is non-polynomial, so that $\bm{u}(\bm{x},0) = \sin(2\pi x)\cos(2\pi y)$ together with periodic boundary conditions on all edges, so that the solution propagates indefinitely. Regular grids are constructed using between $81$ and \numprint{22500} quadrilateral elements in the conformal case. The non-conformal case incorporates two interfaces to ensure that the periodic boundaries are conformal to one another for ease of implementation. This results in three sub-domains, with the central one shifted vertically by half a cell length to create a non-conformal grid. An example non-conformal grid is shown in figure~\ref{fig:vortexMeshNonConformal}. For a given value of $h$, the non-conformal mesh will therefore have a slightly higher number of elements than its conformal counterpart. Polynomial orders of $P=3$ through $P=11$ are considered for each grid, and we select $Q=P+2$ quadrature points in each coordinate direction so as to exactly integrate the mass matrix and remove any spurious aliasing error due to the use of numerical integration. We select a timestep size of $\Delta t = 10^{-3}$ and measure the error after one tenth of a cycle (i.e. $t = 1$) so that error due to timestepping is reduced. In addition, for the point-to-point interpolation method we also investigate an additional setup with $Q=2P+2$ number of quadrature points to investigate the effect of dealiasing on the resulting error and to mirror the strategies employed in later sections.

\begin{table}
\sisetup{round-mode=places, round-precision=2, table-format=2.2}
\centering
\caption{Convergence rates for the conformal case, the non-conformal mortar method, and the non-conformal point-to-point interpolation method.}
\label{tab:convergence}
\begin{tabular}{@{}lSSSS@{}}
\toprule
               			& \multicolumn{4}{c}{Convergence rates for order $p$}                            \\ \cmidrule(l){2-5} 
Poly. order     	    & {$3$} 		        & {$5$}		        & {$7$}		        & {$9$}		     \\ \midrule
Conformal       		& 5.190344294791632 & 7.026373683930996	& 8.878617042670665	& 10.749923787852943 \\
Mortar             		& 4.690484678864176	& 6.302550568111699	& 8.2603084728647	& 9.98215522400944	 \\
Interp, $Q=P+2$         & 4.085677289494459	& 5.885594721109839 & 7.781001409045954 & 9.561802452077583  \\ 
Interp, $Q=2P+2$        & 4.33320093456086	& 5.995077846467883 & 8.054200150323062 & 9.978753147105634  \\ \bottomrule

\end{tabular}
\end{table}

Figure~\ref{fig:convergence} and Table~\ref{tab:convergence} highlight the convergence properties in the $L^2$ sense of the two non-conformal methods, together with the conformal interface. To ensure clarity the results in the graph have been trimmed to remove points from each polynomial order after the minimum $L^2$ error has been reached owing to the finite precision being used. Convergence rates are approximated from the gradients of curves in Figure~\ref{fig:convergence}, asides from at $p=11$ where this is omitted due to lack of data points. The results of this study show that for smooth solutions and at higher orders the mortar method, dealiased interpolation method and conformal cases all show near identical results. These results ratify both that the solvers are implemented correctly and, moreover, that both non-conformal interface handling methods yield similar convergence rates of around $P+1$, while in the conformal setting the rates are $P+2$.

\subsection{Decay properties}

\begin{figure}
    \centering
        \begin{subfigure}[b]{0.98\textwidth}
        \centering
        \includegraphics[width=0.6\textwidth]{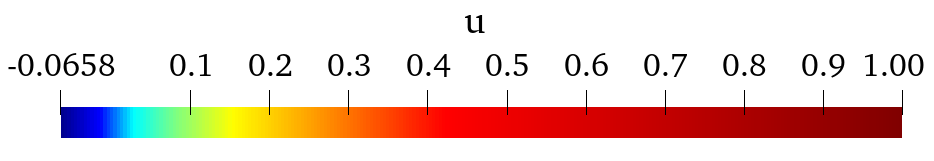}
        \label{fig:gaussianScale}
    \end{subfigure}
    \begin{subfigure}[b]{0.49\textwidth}
        \centering
        \includegraphics[width=1.0\textwidth]{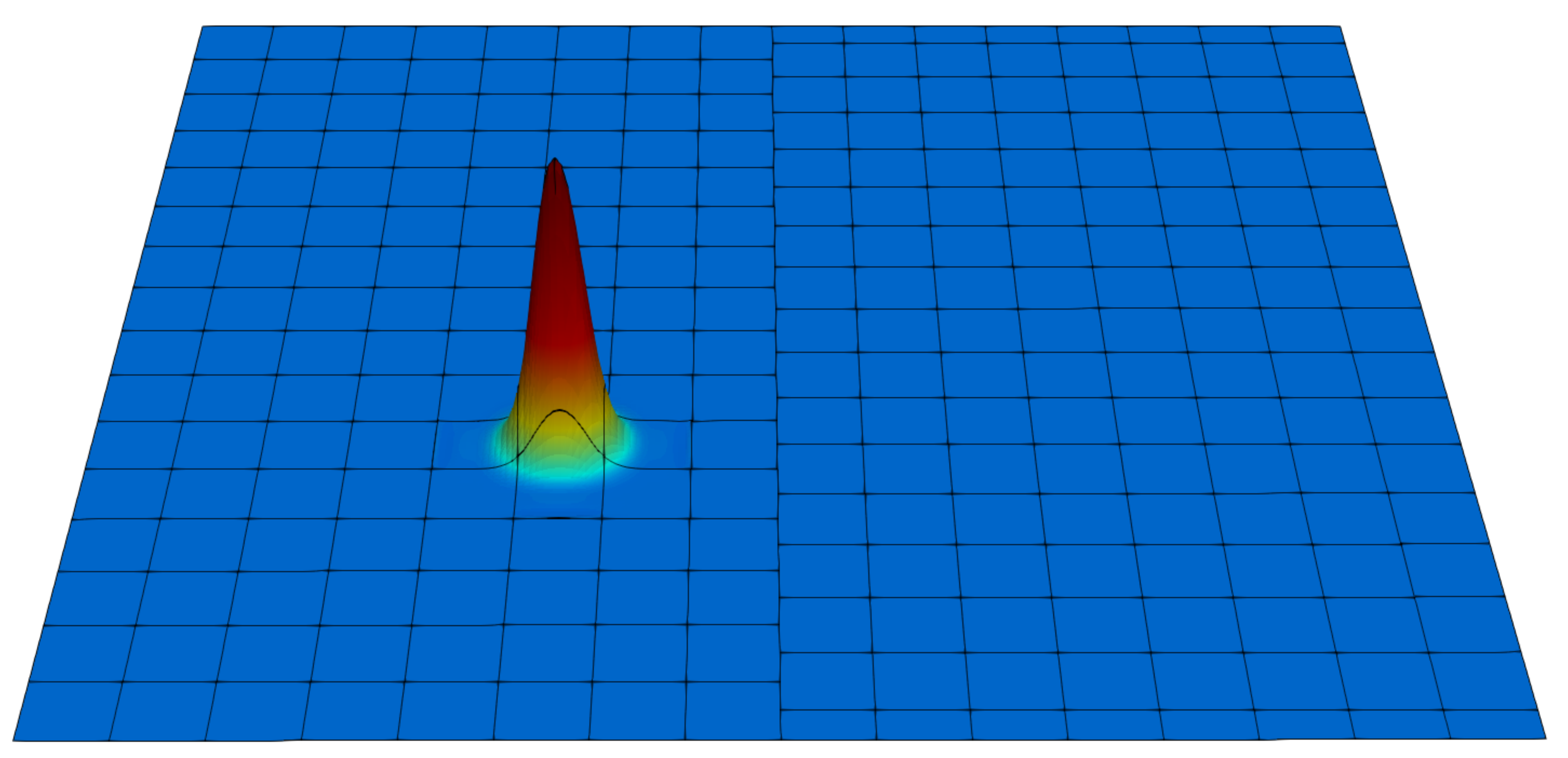}
        \caption{at 0 cycles ($t = 0$)}
        \label{fig:gaussian0}
    \end{subfigure}
    \begin{subfigure}[b]{0.49\textwidth}
        \centering
        \includegraphics[width=1.0\textwidth]{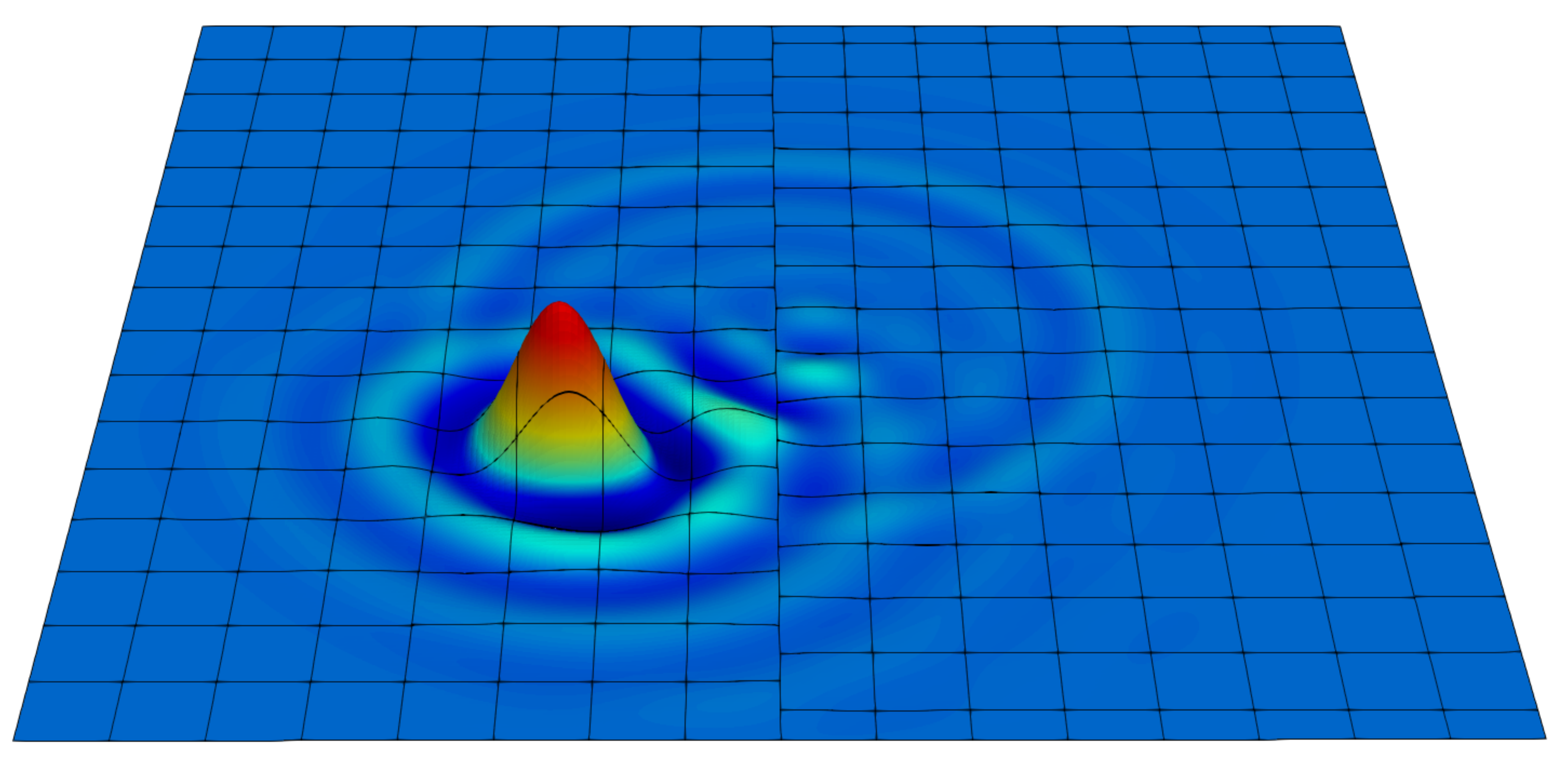}
        \caption{at 100 cycles ($t = 200\pi$)}
        \label{fig:gaussian100}
    \end{subfigure}
    \caption{The decay in the Gaussian peak over the 100 cycles for the 4\textsuperscript{th} order point-to-point interpolation method.}
    \label{fig:gaussMesh}
\end{figure}

In order to more robustly validate each method, we now consider a more challenging problem at varying degrees of resolution. In order to evaluate the numerical diffusion that is introduced by the presence of an interface, we consider the rotation of a Gaussian in a circular manner using the transport equation. More precisely we utilise the same transport equation as the previous setting but now consider the velocity $\bm{v}(x,y) = (-x,y)$, so that the initial scalar Gaussian field $u(\bm{x})_{t=0} = \exp(-\left\|\bm{x}-\bm{x}_0\right\|^2/\sigma^2)$ is rotated around the origin. We consider a domain $\Omega=[-2,2]$, using an initial starting point $\bm{x}_0 = (-0.625, -0.625)$ with $\sigma=0.1$. The mesh used in this test consists of $16 \times 16$ uniformly-sized quadrilateral cells for the conformal case, while in the non-conformal cases the right-hand half of the grid has been displaced by half a cell vertically along the central interface in relation to the left-hand side, as shown in figure~\ref{fig:gaussMesh}. Constructing the mesh in this way aims to keep a consistent cell density by ensuring the half cell height sections are on the extreme ends of the interface, distant from where the peak crosses the interface. We also note that the selection of $\bm{x}_0$ is designed to place the peak in the centre of a cell to the left of the interface, this is to ensure minimal interaction with the domain boundaries which all have a homogeneous Dirichlet condition imposed on them.

The peak starts at $t = 0$ with a maximum value of 1. Unlike the exact solution, which precisely preserves this peak indefinitely, the non-polynomial nature of the solution field means that we can expect the peak to decrease every rotational cycle due to numerical diffusion introduced by each method. We then measure the $L^\infty$ norm of the solution precisely through a minimisation problem -- i.e. we do not solely sample the error at quadrature points, as at lower orders very few quadrature points are used within each element, and this may lead to a significant difference in the observed error. We select a timestep size of $\Delta t = 10^{-3}$, and for each combination of polynomial order and interface handling method, we measure the $L^\infty$ norm after 100 cycles of the Gaussian (i.e. $t=200\pi$). We note that at lower polynomial orders, the solution will be underresolved by design -- our aim in this series of simulations is to examine how this affects numerical stability across the methods and/or if there are significant differences in performance of the methods. The results of these experiments are shown in table~\ref{tab:gaussianResults}.

\begin{table}[h]
\sisetup{round-mode=places, round-precision=4, table-format=1.4}
\centering
\caption{Results of Gaussian peak value after 100 cycles for varying basis orders}
\label{tab:gaussianResults}
\begin{tabular}{@{}lSSSSSSS@{}}
\toprule
               			& \multicolumn{7}{c}{$\| \cdot \|_{L^{\infty}}$ after 100 cycles} \\ \cmidrule(l){2-8} 
Poly. order     	    & {4} 		& {5}		& {6}		& {7}		& {8}		& {9} 		& {10} 		\\ \midrule
Conformal       		& 0.440289 	& 0.687013	& 0.829138	& 0.930197	& 0.975797	& 0.993267	& 0.998420	\\
Mortar             		& 0.438670	& 0.681912	& 0.828588	& 0.929262	& 0.975604	& 0.993124	& 0.998434	\\
Point-to-point  		& 0.481175	& 0.705565 	& 0.839422  & 0.947575  & 0.972019  & 0.993560	& 0.998014	\\ \bottomrule

\end{tabular}
\end{table}

Two trends are evident from the presented results. At the lowest polynomial order of $P=4$, we see a reasonable level of difference in the point-to-point method vs. the mortar and conformal grids. The oscillations of numerical error at this order are clearly evident in figure~\ref{fig:gaussian100} which shows the 4\textsuperscript{th} order point-to-point interpolation method after 100 cycles. Curiously, the values observed at $P=4$ through $P=7$ are higher for the point-to-point method than both the mortar/conformal cases, indicating that the point-to-point method is somewhat better able to resolve the peak of the Gaussian. It is also clear that as the  polynomial order increases this difference in the maximum value decreases, so that at $P\geq 8$ the results are essentially identical. Broadly speaking, however, the performance of the methods is reasonably comparable across the range of polynomial orders.

\subsection{Long-time advection of an isentropic vortex}

In order to examine the non-conformal methods in more realistic problems, whilst still considering their long-term stability and diffusion properties, we now move on to consider a nonlinear problem. In particular, we consider the compressible Euler equations in two dimensions. In this case, the conserved variables are given as $\bm{u} = [\rho, \rho u, \rho v, E]$ with $\rho$ being the density, $(u,v)$ the fluid velocity, $E$ the specific total energy, and
\begin{align*}
\bm{F}(\bm{u}) = \left[
  \begin{array}{cc}
    \rho u       & \rho v       \\
    p + \rho u^2 & \rho uv      \\
    \rho uv      & \rho v^2 + p \\ 
    u (E + p)    & u(E + p)
  \end{array}
\right],
\end{align*}
where $p$ is the pressure. To close the system we need to specify an equation of state; in this case we use the ideal gas law $p=\rho RT$ where $T$ is the temperature and $R$ is the gas constant.

To consider long-term stability, we opt to study an isentropic vortex that is advected at constant velocity through periodic boundaries. This is a commonly used benchmark when testing numerical discretisation of the compressible Euler equations, particularly for higher-order codes, as it is one of the few problems that admits an exact solution calculable at all times whilst also being relatively simple to implement~\cite{deGrazia2014, wang2007, nasa1, nasa2}.

\begin{figure}
    \centering
    \begin{subfigure}[b]{0.495\textwidth}
        \centering
        \includegraphics[width=\textwidth]{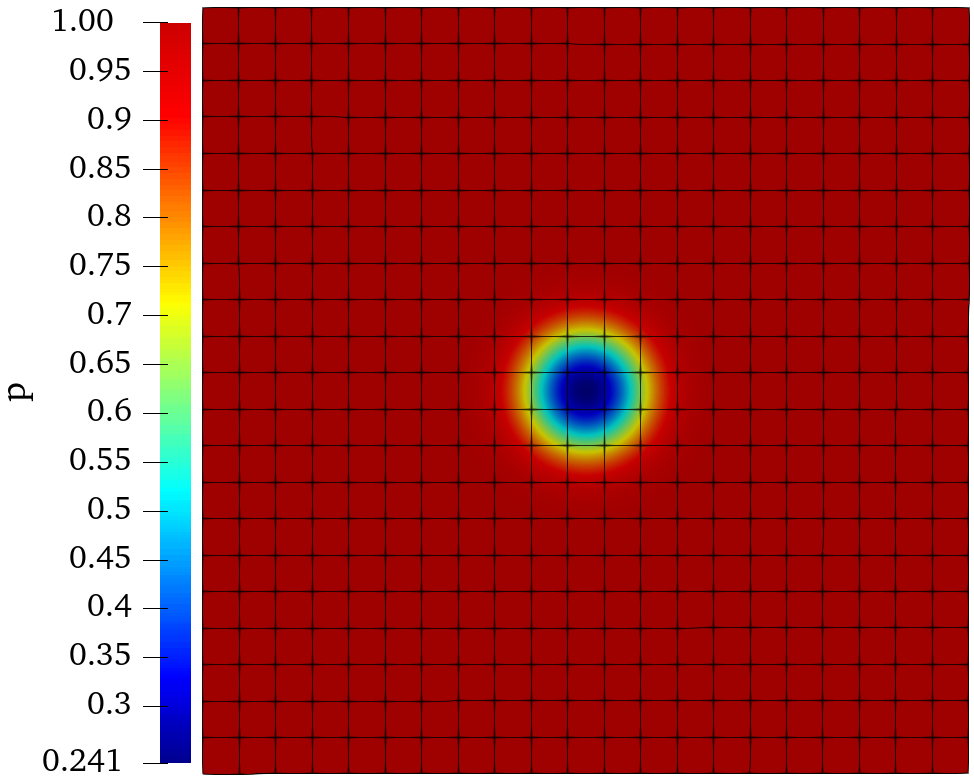}
        \caption{Conformal mesh }
        \label{fig:vortexMeshConformal}
    \end{subfigure}
    \begin{subfigure}[b]{0.495\textwidth}
        \centering
        \includegraphics[width=0.797\textwidth]{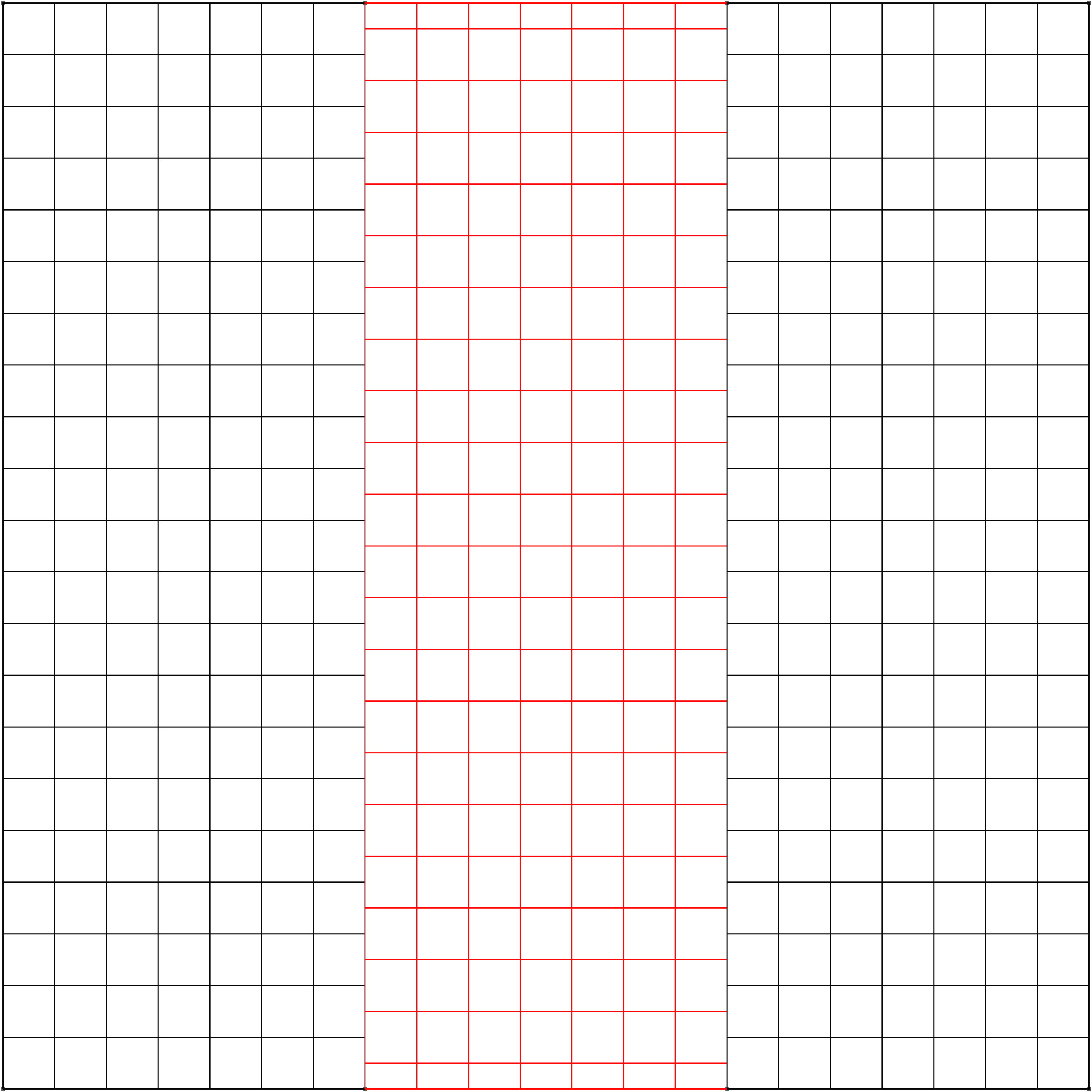}
        \caption{Non-conformal mesh}
        \label{fig:vortexMeshNonConformal}
    \end{subfigure}
    \caption{The two meshes used for the isentropic Euler vortex cases, (a) the conformal case has the initial projection of the pressure field overlaid.}
    \label{fig:vortexMesh}
\end{figure}

For our purposes, we consider a domain $\Omega = [-5,5]^2$. At any given time $t$, the solution for the isentropic vortex is given by the equations
\begin{align}
\begin{split}
    \rho &= \left(1 - \frac{\beta^{2}\left(\gamma - 1\right)\text{e}^{2f}}{16\gamma\pi^{2}}\right)^{\frac{1}{\gamma-1}}, \\
    u &= \left(u_{0} - \frac{\beta \text{e}^{f}\left(y-y_{0}\right)}{2\pi}\right), \\ 
    v &= \left(v_{0} - \frac{\beta \text{e}^{f}\left(x-x_{0}\right)}{2\pi}\right), \\ 
    E &= \frac{\rho^{\gamma}}{\gamma - 1} + \frac{1}{2}\rho\left(u^{2} + v^{2}\right),
\end{split}
\label{eqn:vortex}
\end{align}
where $f = 1 - \left(x - x_{0}\right)^{2} + \left(y - y_{0}\right)^{2}$. We select an initial vortex position $(x_0, y_0) = (0,0)$ with strength $\beta=5$ and $\gamma=1.4$, and advect the velocity in the $x$-direction with velocity $(u_0, v_0) = (1,0)$. The initial vortex size and location can be seen in figure~\ref{fig:vortexMeshConformal}.

The concept behind this series of simulations is much the same as the rotating Gaussian peak; i.e. we wish to cycle the vortex through the domain a number of times, and compare the error as a function of time for each interface method. This will be undertaken for polynomial orders, $P=3$ through $P=7$ on a grid of fixed size, where the lower orders are expected to be under-resolved and the higher orders somewhat more resolved. To impose this, a single pair of periodic conditions at the constant $x$ boundaries were used so that $\bm{u}(-5,y)=\bm{u}(5,y)$, while the constant $y$ boundaries were set to free-stream conditions. Although it may seem more natural to impose periodic boundaries on all of the edges of the domain, in a similar fashion to~\cite{nasa1} we found that this leads to a gradual accumulation of numerical error, which left unchecked eventually causes simulations to diverge. The solution, proposed in~\cite{nasa1} and~\cite{nasa2}, is to impose farfield conditions at constant $y$ boundaries, which allows recirculated waves of accumulated numerical error to escape the domain and avoid premature divergence. This is particularly important in this case, as in order to further reduce sources of artificial dissipation, we  elect to use the exact Riemann solver of Toro~\cite{toro-2009} to calculate the numerical flux $\tilde{\bm{f}}(\bm{u}^+, \bm{u}^-)$. We note that although this is computationally expensive, cheaper solvers such as the Roe solver may introduce additional numerical diffusion \cite{li2017}.

One additional consideration that needs to be taken in this nonlinear regime is the order of integration used to evaluate integrals in the weak form of eq.~\eqref{eq:weak}. Aliasing errors are a well known phenomenon in this regime, due to the cubic nonlinearity that arises in the definition of the Euler equations, as well as the non-polynomial flux term which calculated between elements~\cite{deGrazia2014}. Typically it is necessary to use a higher order of quadrature than is used for linear problems, in order to remove sources of aliasing error due to under-integration of these terms. For this reason, we consider two different numbers of quadrature points with $Q = P + 2$ and $Q = 2P + 2$, respectively.

The conformal case consists of a singular domain made up of a $21 \times 21$ regular quadrilateral mesh as shown in figure~\ref{fig:vortexMeshConformal}. The resulting non-conformal mesh consists of the domains $7 \times 21$, $7 \times 22$ and $7 \times 21$ elements, as shown in figure~\ref{fig:vortexMeshNonConformal}.
The periodic boundary also allows us to conveniently express the time in cycles, where one cycle is the length of time taken for the vortex to propagate through the domain and return to its initial position. In our case, with a propagation speed of $u_0=1$ and domain length $L=10$, this leads to the same exact solution every $t=10$. We can calculate the exact solution at any time, $t$, by moving the vortex centre by $u_{0}t$ in the $x$-direction and making sure to account for the periodic condition. We select a fixed timestep size of $\Delta t=10^{-3}$ and use the same explicit 4th-order Runge-Kutta timestepping scheme as in previous results.

\begin{figure}
    \centering
    \includegraphics[width=\textwidth]{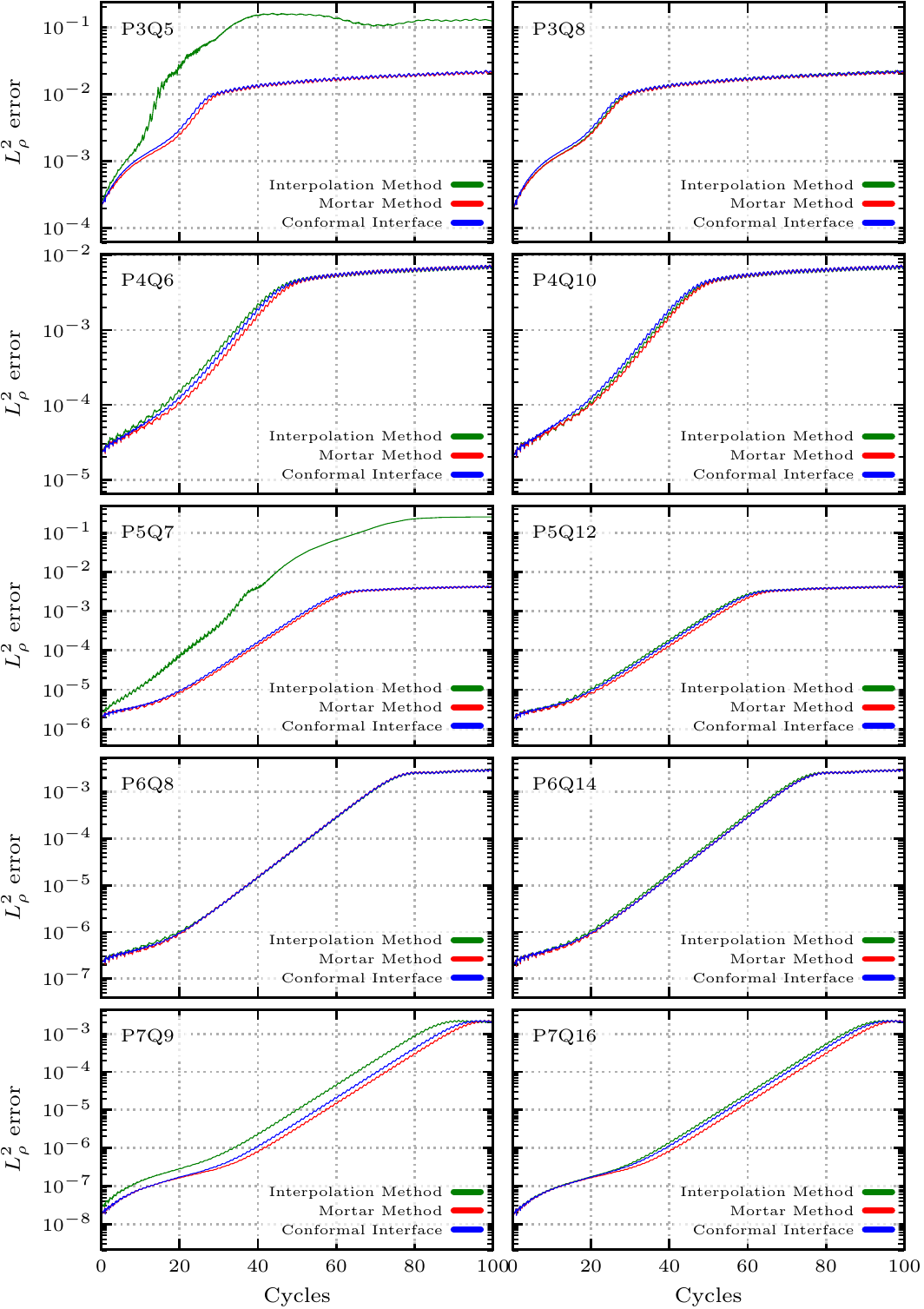}
    \caption{The evolution of the $L_{\rho}^{2}$ error as a function of cycles of the vortex through the domain. Rows of the figure denote polynomial order $P=3$ through $P=7$, and columns denote integration orders $Q=P+2$ and $Q=2P+2$, respectively.}
    \label{fig:vortexErrors}
\end{figure}

In figure~\ref{fig:vortexErrors}, we visualise the $L^2$ error of the density field $\rho$, denoted by $L_{\rho}^{2}$, for a simulation spanning 100 cycles of the vortex through the domain. This figure yields a number of interesting features that warrant further discussion. Firstly, as validation of our results, we note that the broad characteristics of the conformal error broadly agree with those seen in other work and, in particular, those of~\cite{nasa1}. More generally, we observe that the conformal method and mortar method yield extremely close results for all polynomial and quadrature orders under observation, which we would expect given the similar levels of resolution and the local conservation properties of the mortar method.

However, when considering the point-to-point interpolation method, there are indeed clear differences in comparison to the mortar and conformal methods. Perhaps the most obvious peculiarity is the possible relationship between odd numbers of quadrature points and the long term stability of the interpolation method; P3Q5 and P5Q7 in figure~\ref{fig:vortexErrors} show significant divergence at low cycle counts. The errors in this case appear to be related to aliasing error: as the integration order is increased to $Q=8$ and $Q=12$ respectively, the results remain consistent with those found by the mortar method and the reference conformal case.

\begin{figure}
    \centering
    \includegraphics[width=0.94\textwidth]{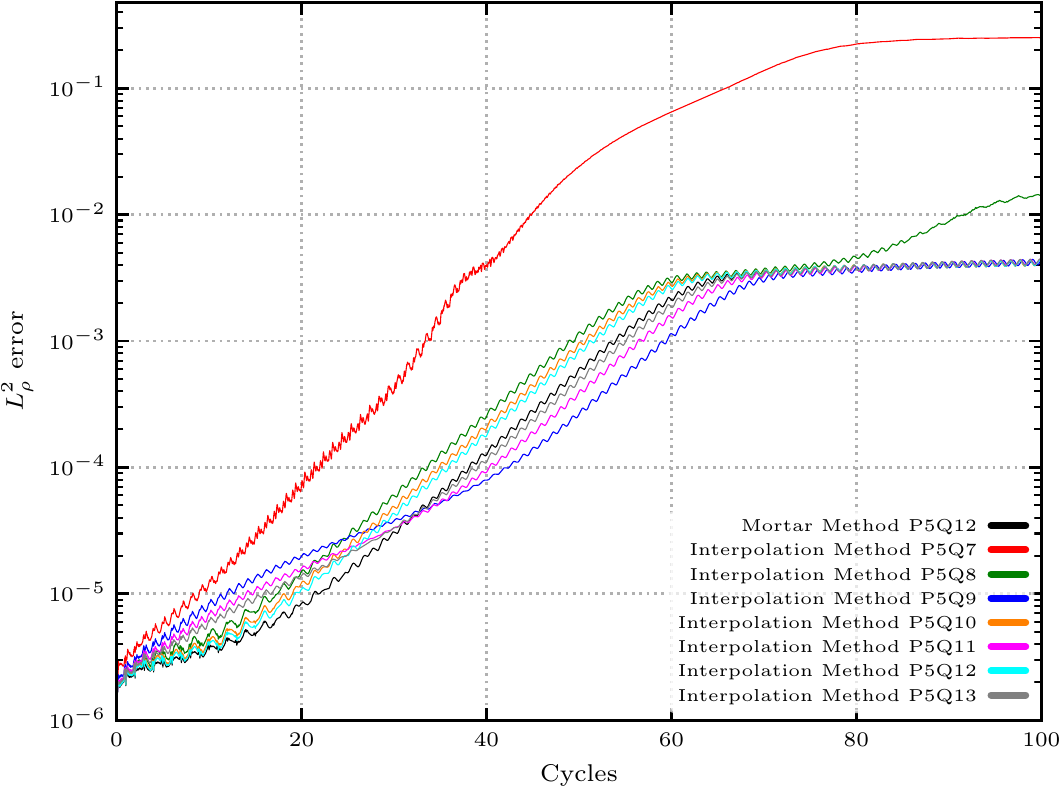}
    \caption{The evolution of the $L_{\rho}^{2}$ error over time for $P=5$ with the interpolation method at various quadrature point numbers. P5Q12 for the mortar method is also shown as a comparison.}
    \label{fig:vortexErrorsP5}
\end{figure}

\begin{figure}
    \centering
    \includegraphics[width=0.94\textwidth]{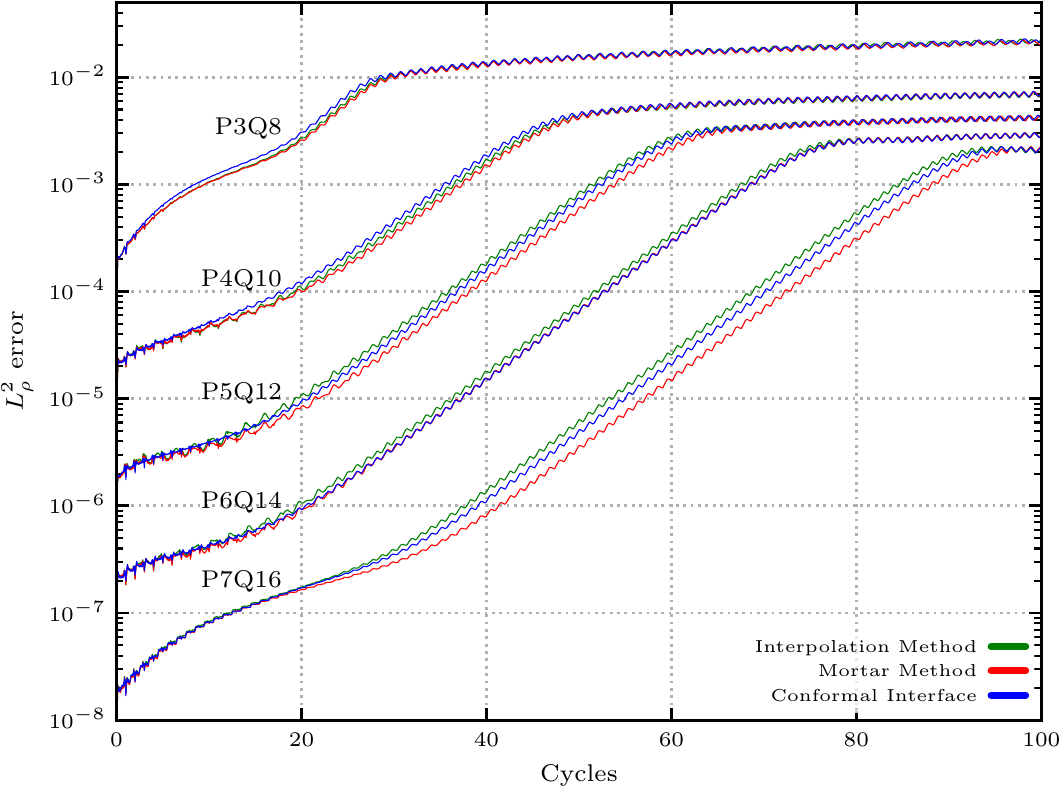}
    \caption{The evolution of the $L_{\rho}^{2}$ error over time for P3, P4, P5, P6 and P7 with $Q=2P + 2$.}
    \label{fig:vortexErrorsOne}
\end{figure}

To investigate the effect of aliasing and integration order further, additional point-to-point interpolation simulations were run at $P=5$ with quadrature orders ranging between $Q=7$ and $Q=13$. Figure~\ref{fig:vortexErrorsP5} depicts the $L_\rho^2$ error for these cases. The pronounced abnormality at $Q=7$ is clearly visible, and indeed at $Q=8$, there is a sudden increase in error after $\sim 70$ cycles which is indicative of further long-time increases in error. However, for $Q\geq 9$, we observe much more consistent trends and better agreement with the mortar and conformal cases. More generally then, we can state that so long as appropriate levels of aliasing are used so that $Q=2P+2$, the interpolation method closely follows the same trend as the mortar method and the benchmark conformal case, with the same reduction in error as polynomial order increases. These cases are all visualised in a single figure~\ref{fig:vortexErrorsOne} to highlight this more clearly.


\begin{table}[h]
\sisetup{round-mode=places, round-precision=4, table-format=1.4}
\centering
\caption{Computational costs for the isentropic Euler vortex with $Q=P+2$ cases.}
\label{tab:compCosts}
\begin{tabular}{@{}lSSSSS@{}}
\toprule
               	& \multicolumn{5}{c}{Average cost per timestep (s)} \\ \cmidrule(l){2-6} 
Case	      	& {P3Q5} 	& {P4Q6}	& {P5Q7}	& {P6Q8}	& {P7Q9}	\\ \midrule
Conformal       & 0.0201855	& 0.0298848	& 0.0429752 & 0.0577017 & 0.079429  \\
Point-to-point  & 0.0224373	& 0.0325178 & 0.0463944 & 0.0617084 & 0.0679356 \\ 
Mortar          & 0.0872229	& 0.098268	& 0.113477	& 0.129721  & 0.17081   \\ \bottomrule 
\end{tabular}
\end{table}

\begin{table}[h]
\sisetup{round-mode=places, round-precision=4, table-format=1.4}
\centering
\caption{Computational costs for the isentropic Euler vortex with $Q=2P+2$ cases.}
\label{tab:compCosts2}
\begin{tabular}{@{}lSSSSS@{}}
\toprule
               	& \multicolumn{5}{c}{Average cost per timestep (s)} \\ \cmidrule(l){2-6} 
Case	      	& {P3Q8} 	& {P4Q10}	& {P5Q12}	& {P6Q14}   & {P7Q16}	\\ \midrule
Conformal       & 0.035941  & 0.0555044 & 0.0847885 & 0.117629  & 0.164616  \\
Point-to-point  & 0.0366097 & 0.0558076 & 0.0848238 & 0.118458  & 0.163842	\\ 
Mortar          & 0.111157  & 0.130747  & 0.187239  & 0.232467	& 0.269642	\\ \bottomrule
\end{tabular}
\end{table}

In addition we also used this problem to investigate the computational costs associated with each interface method. These simulations were run on a single core of a dual-socket Intel Xeon Gold 5120 system, equipped with 256GB of RAM, with the solver pinned to a specific core in order to reduce the influence of kernel core and socket reassignment mid-process. The time taken per timestep for the $Q=P+2$ cases is shown in table~\ref{tab:compCosts}, and the $Q=2P+2$ cases in table~\ref{tab:compCosts2}. This shows that the conformal and point-to-point interpolation simulation timings are very similar with a small percentage cost associated with the interpolation. Of note are the $P=7$ results, which show the point-to-point interpolation cost as lower than the equivalent conformal case. A further investigation of this case showed that the small variation in the number of elements leads to a reduction in walltime of the evaluation of DG volume terms, possibly owing to the strategies used in Nektar++ to evaluate collective operations more effectively~\cite{moxey2016}. The mortar method shows a larger cost, which decreases in proportion to the other methods as the polynomial order and quadrature point number is increased. For example, at P3Q5 it is over four times as expensive, whilst at P7Q16 it is less than double the conformal cost. This suggests that in the mortar method cases, at least in this setup, the total computational cost is dominated by the interface handling. This is to be expected, since the projection both to and from mortars involves more costly operations than a straightforward interpolation in the point-to-point approach. We note here that although we have taken some steps to optimise the mortar method implementation (e.g. through the caching of $\mathbf{S}$ matrices defined in the previous section), more in-depth techniques such as those found in~\cite{durrwachter2021} may help in reducing walltime for parallel execution in particular.

\section{Extension to three-dimensional cases}
\label{sec:3d}

In this section we consider the extension of our two-dimensional simulations to a realistic three-dimensional fluid dynamics case. In particular, we consider the extension to the compressible Navier-Stokes equations, which in conservative form may be written as
\begin{equation}
\diffp{\bm{u}}{t} + \nabla\cdot\bm{F}(\bm{u}) =
\nabla\cdot\bm{F}_v(\bm{u}, \nabla\bm{u}), \label{eq:ns}
\end{equation}
where $\bm{u} = [\rho, \rho u, \rho v, \rho w, E]$ is the vector of conserved variables in terms of density $\rho$, velocity $\bm{v} = (u_1,u_2,u_3) = (u,v,w)$ and $E$ is the specific total energy. In three dimensions we have that
\begin{align*}
\bm{F}(\bm{u}) = \left[
  \begin{array}{ccc}
    \rho u       & \rho v       & \rho w \\
    p + \rho u^2 & \rho uv      & \rho uw \\
    \rho uv      & \rho v^2 + p & \rho vw \\ 
    \rho uw      & \rho vw      & \rho w^2 + p \\
    u (E + p)    & u(E + p)     & v(E + p)
  \end{array}
\right],
\end{align*}
We again use the ideal gas to close the system. The tensor of viscous
forces $\bm{F}_v(\bm{u}, \nabla\bm{u})$ is defined as
\begin{align*}
  \bm{F}_v(\bm{u}, \nabla\bm{u}) &= \left[
    \begin{array}{ccc}
      0 & 0 & 0\\
      \tau_{xx} & \tau_{yx} & \tau_{zx} \\
      \tau_{xy} & \tau_{yy} & \tau_{zy} \\
      \tau_{xz} & \tau_{yz} & \tau_{zz} \\
      A & B & C
    \end{array}\right],
\end{align*}
with
\begin{align*}
  A &=  u\tau_{xx} + v\tau_{xy} + w\tau_{xz} + k\partial_x T,\\
  B &=  u\tau_{yx} + v\tau_{yy} + w\tau_{yz} + k\partial_y T,\\
  C &=  u\tau_{zx} + v\tau_{zy} + w\tau_{zz} + k\partial_z T,
\end{align*}
where in tensor notation the stress tensor $\tau_{x_ix_j} = 2\mu (\partial_{x_i}u_i + \partial_{x_i}u_j - \tfrac{1}{3}\partial_{x_k}u_k\delta_{ij})$, $\mu$ is the dynamic viscosity calculated using Sutherland's law, $k$ is the thermal conductivity and $\delta_{ij}$ is the Kronecker delta.

From a numerical perspective, we adopt the same discontinuous Galerkin formulation to discretise equation~\eqref{eq:ns}. However we note that the inclusion of the viscous term $\bm{F}_v(\bm{u}, \nabla\bm{u})$ requires additional treatment, in particular a careful selection of flux terms in order to preserve spatial accuracy. In the simulations below, we adopt the local discontinuous Galerkin (LDG) approach, wherein an auxiliary variable $\bm{q} = \nabla\bm{u}$ is introduced and discretised alongside equation~\eqref{eq:ns}. With a careful choice of alternating fluxes (so that $\tilde{\bm{q}} = \bm{q}^+$ and $\tilde{\bm{u}} = \bm{u}^-$, or vice versa), high-order accuracy can be preserved~\cite{cockburn1998}. 

\subsection{Implementation considerations}

In order to extend the formulation in section~\ref{sec:theory} to three dimensions for a non-conformal Navier-Stokes simulation, consideration has to be given to a number of implementation changes, which we briefly outline in this section.

Both the mortar and point-to-point interpolation method require the evaluation of the solution at arbitrary points within the skeleton of the mesh. In 2D simulations, this requires evaluation within an interval; however in 3D this could conceivably be evaluated in either quadrilateral or triangular faces, depending on the element type: for example hexahedra possess purely quadrilateral faces, whereas tetrahedra possess triangular faces. Although the barycentric evaluation approach proposed in section~\ref{sec:theory} naturally extends to higher-dimensional quadrilaterals and hexahedra through a tensor product of one-dimensional evaluations, for triangular elements and other three-dimensional shapes, most discontinuous Galerkin implementations documented in the literature are based around the selection of a set of cubature points combined with Lagrange interpolants as basis functions. Typical examples of such distributions are the Fekete~\cite{taylor2000} or electrostatic points~\cite{hesthaven1998}, which both provide better conditioning of operators when compared to evenly-spaced points~\cite{karniadakis2005}. In this case, interpolation may be done via the typical route of generating an interpolation matrix with the aide of a Vandermonde matrix, as outlined in~\cite{hesthaven2008}.

However, we note that in the spectral/$hp$ formulation of Karniadakis \& Sherwin~\cite{karniadakis2005} which forms the numerical basis for Nektar++, higher $d$-dimensional simplicies, as well as other hybrid shape types such as prisms and tetrahedra, are represented instead on a collapsed coordinate space, denoted by $\bm{\eta}\in[-1,1]^d$. We give a brief overview of this formulation here, leaving further details to the aforementioned reference. Each collapsed coordinate spaces is mapped to the desired reference elemental shape through the use of Duffy transformations. For example, in a triangular element with reference coordinates $\{ (\xi_1, \xi_2)\ |\ \xi_1,\xi_2 \geq -1, \xi_1+\xi_2 \leq 0 \}$, we have that
\[ 
\eta_1 = 2\frac{1+\xi_1}{1-\xi_2}-1, \quad \eta_2 = \xi_2
\]
Evaluation of quadrature therefore occurs on the collapsed space which, being a quadrilateral or hexahedron in two or three dimensions, may be equipped with a set of tensor-product integration points. Typically, this is chosen to be a set of Gauss-Lobatto points in the $\eta_1$ direction, and Gauss-Radau points in the $\eta_2$ direction in order to avoid explicit evaluations near the (removable) singularity in the Duffy transformation which occurs when $\xi_2 = 1$. In this manner, barycentric interpolation can still be applied in order to increase the computational speed of the simulation for triangular elements, as well as other higher-dimensional shapes. This is a topic of broader interest and under investigation in a wider range of areas~\cite{laughton2021}.

Another significant issue to overcome in the case of mortaring is the construction of the mortar space. As noted in the preceding sections, for generic interfaces between unstructured grids, this can pose a significant challenge, although several techniques have been demonstrated in the literature to handle cases where the geometry is extruded and thus elicits structure that can be exploited. In the remainder of this section, we opt therefore to consider only the point-to-point interpolation approach, since the preceding section clearly demonstrates near-identical behaviour when compared to normal conformal simulations, and the aforementioned references demonstrate the viability of this approach in three-dimensional simulations. The more pertinent question is therefore how the point-to-point interpolation approach performs in this setting, which has yet to be examined in these cases to the best of the authors' knowledge.

\subsection{Simulation of a Taylor-Green vortex}

In order to examine the performance of the point-to-point interpolation method, we consider the simulation of a Taylor-Green vortex at a Reynolds number $\mathrm{Re}=\numprint{1600}$, which has become a benchmark case for the evaluation of higher-order CFD codes. In this case, starting vortices are defined in a periodic box $\Omega = [-L\pi, L\pi]^3$, given a reference length $L$, which break down into turbulent eddies before decaying due to viscous effects. The initial conditions are given in primitive variables $(\bm{v}, p)$ as
\begin{align*}
    u &= V_\infty \sin(x/L)\cos(y/L)\cos(z/L), \\
    v &= -V_\infty \cos(x/L)\sin(y/L)\cos(z/L), \\
    w &= 0, \\
    p &= \rho_\infty V_\infty^2\left[ \frac{1}{\gamma \mathrm{Ma}_\infty^2} + \frac{1}{16}\left(\cos(2x/L) + \cos(2y/L)\right)\cdot\left(\cos(2z/L) + 2\right) \right],
\end{align*}
with the Reynolds number $\mathrm{Re} = \rho_\infty U_\infty L/\mu$ and the Prandtl number $\mathrm{Pr} = 0.71$. Although the Taylor-Green vortex problem is traditionally examined in the setting of an incompressible flow, we approach this limit by considering flows with low compressibility effects so that the Mach number $\mathrm{Ma}_\infty = 0.1$. A simulation is then conducted across the time interval $t_c\in[0,20]$, where the convective timescale $t_c = tV_\infty/L$. We select an explicit second-order Runga-Kutta time integration scheme, with the timestep adjusted to maintain a Courant-Friedrichs-Lewy (CFL) condition of 0.2.

\subsubsection{Kinetic energy dissipation rate}

A key quantity of interest in this simulation is the evolution of the kinetic energy dissipation rate $\epsilon$, where
\[
\epsilon = -\diff{{E_k}}{t}, \quad E_k = \frac{1}{\rho_\infty|\Omega|} \int_\Omega \tfrac{1}{2}\rho\|\bm{v}\|^2 \,\mathrm{d}\bm{x},
\]
since the peak dissipation is a difficult quantity to resolve closely for under-resolved simulations. Examination of $\epsilon$ therefore gives an indication as to the numerical performance of the scheme and inherent numerical diffusion. For an incompressible fluid, $\epsilon = 2\mu\mathcal{E}/\rho_\infty$ where the enstrophy $\mathcal{E}$ is computed as
\[
\mathcal{E} = \frac{1}{\rho_\infty|\Omega|} \int_\Omega \tfrac{1}{2}\rho\|\bm{\omega}\|^2 \mathrm{d}\bm{x}.
\]
This equality does not strictly hold for a compressible fluid. However, as the additional contributions that appear in the exact expression depend on the divergence of the velocity, for this close-to-incompressible case their contribution is very small and can therefore be omitted.

\begin{figure}[h]
    \centering
    \includegraphics[width=0.9\textwidth]{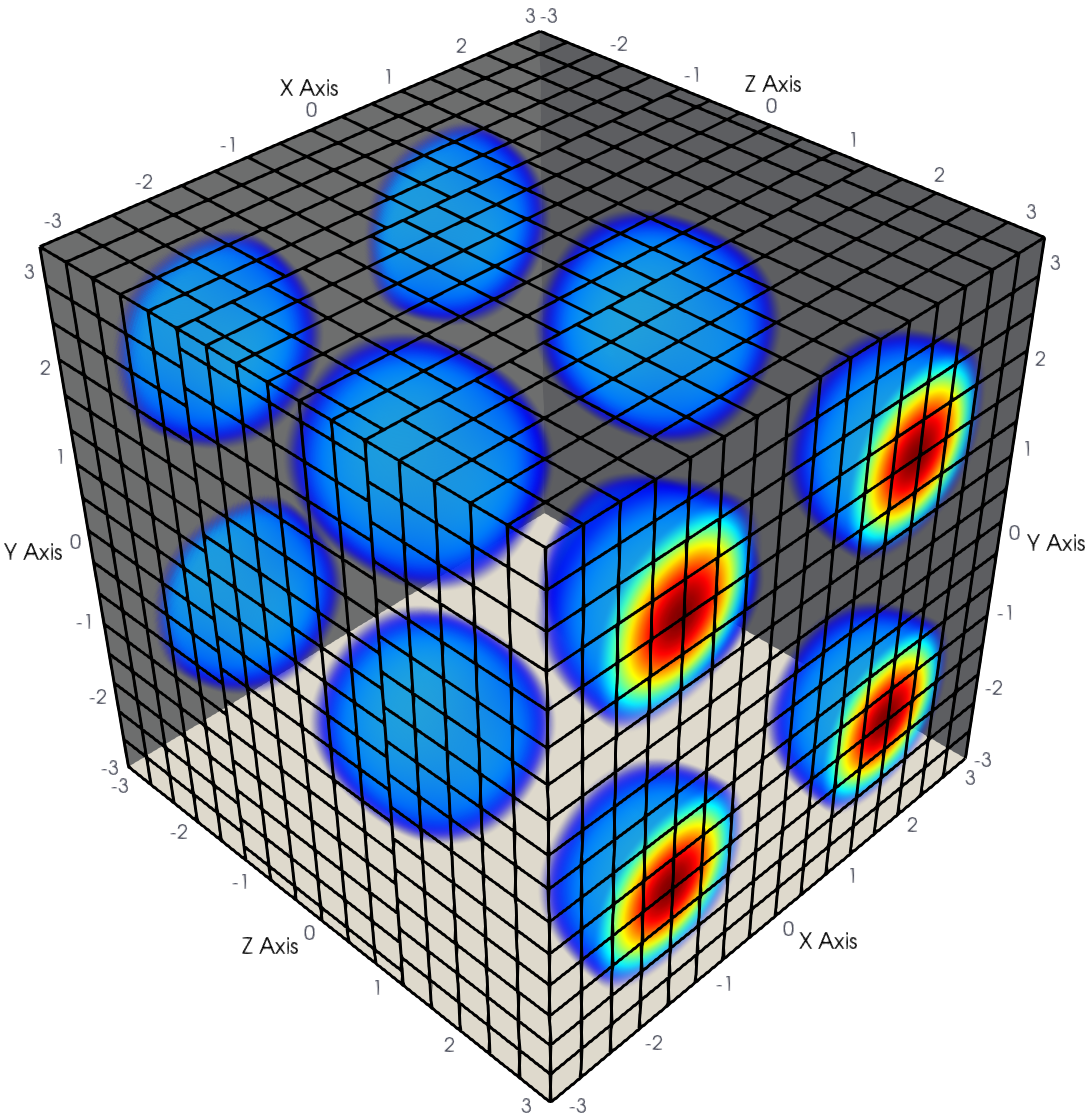}
    \caption{Non-conformal 3D mesh for the $64^3$ DOF Taylor-Green vortex case showing $\|\omega\|$ at $t_c = 0$.}
    \label{fig:tg-mesh}
\end{figure}

\begin{figure}[h]
    \centering
    \includegraphics[width=1.0\textwidth]{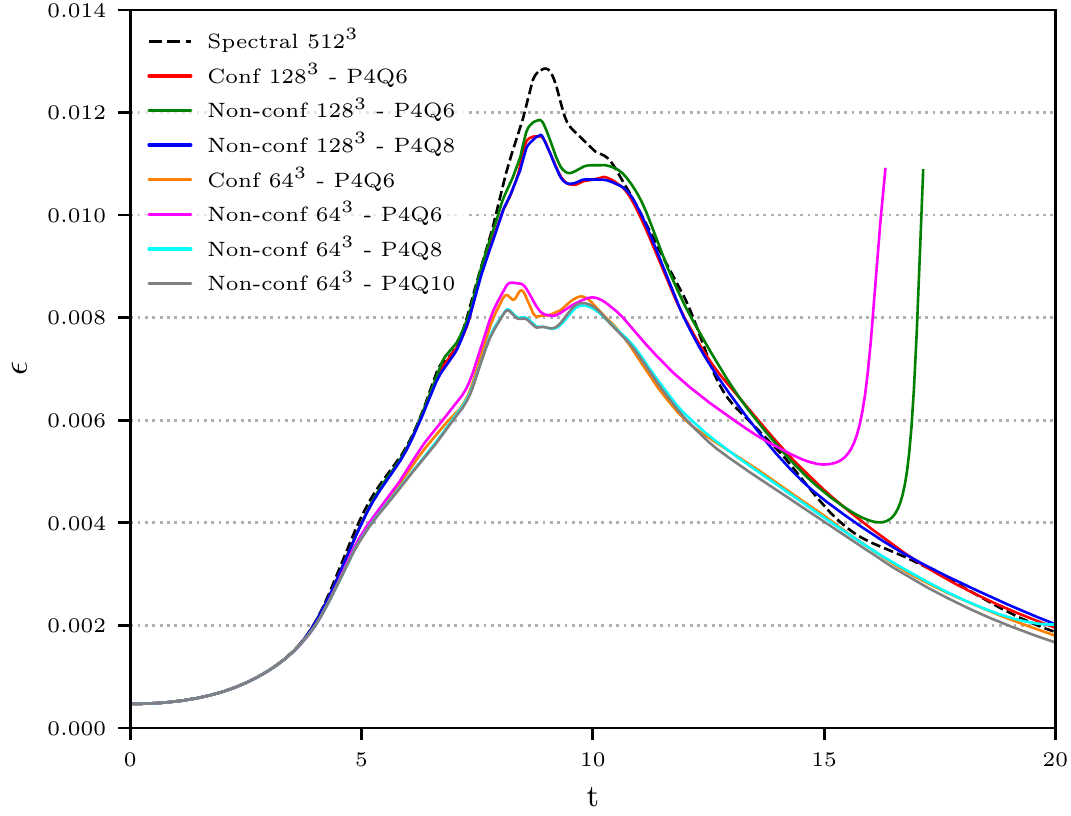}
    \caption{Evolution of the dimensionless energy dissipation rate in the Taylor-Green vortex, showing the reference spectral data and both conformal and non-conformal at $64^3$ and $128^3$ DOF.}
    \label{fig:tg-vortex}
\end{figure}

In figure~\ref{fig:tg-vortex} we visualise the evolution of $\epsilon$ from a number of simulations:
\begin{itemize}
    \item reference DNS data from a spectral simulation using $512^3$ grid points in each spatial direction;
    \item simulations on a conformal mesh with $16^3$ and $32^3$ equally-sized hexahedra at polynomial order $P=4$ and quadrature order $Q=6$, for a total of $64^3$ and $128^3$ degrees of freedom (DOF) equivalent resolution;
    \item simulations on a non-conformal mesh at equivalent levels of resolution and polynomial order to the conformal case, with various levels of dealiasing. The mesh used for the non-conformal simulations as well as the initial vorticity condition is visualised in figure~\ref{fig:tg-mesh}.
\end{itemize}

As is seen from the reference spectral data, the maximum kinetic energy dissipation appears at $t_c\approx 9$. The simulations of the conformal DG case at $128^3$ clearly demonstrate the capability of the discontinuous Galerkin method to broadly resolve all but the very peak of the simulation even at a factor of 4 reduction in resolution, with greater numerical diffusion occurring when the resolution is reduced further to a $64^3$ resolution. These results are broadly in line with a number of other simulations that appear in the literature, for example in~\cite{yan2021, fehn2018, drikakis2007, chapelier2012}.

Of course the central interest of this work is to examine the effect of the point-to-point interpolation in the non-conformal case. Two trends are immediately apparent. Plainly the clearest aspect of  figure~\ref{fig:tg-vortex} is that when run without any polynomial dealiasing, the simulation becomes unstable at $t_c\approx 16$. From a fluid dynamics perspective, this is during the vortex saturation phase, where breakdown of the vortices is approaching the viscous limit and thus the smallest features are starting to appear in the flow. At under-resolution, this highlights the increased fragility of the point-to-point interpolation approach, likely owing to the appearance of oscillatory effects as structures break down and thus regularity of the solution across the interface decreases. However, as we observe in the previous section and figures, the use of moderate levels of dealiasing, commensurate with what is typically leveraged for underresolved simulations, can stabilise the simulation. Moreover, aside from this instability, it is clear that the non-conformal simulations very closely track the evolution of $\epsilon$ when compared to the conformal cases, with only small deviations observed from the conformal cases that can perhaps be attributed to a small difference in number of elements between simulations. Nevertheless, these simulations further emphasise that careful consideration of dealiasing is critical in this setting.

\begin{figure}
        \begin{subfigure}[b]{0.50\textwidth}
        \centering
        \includegraphics[width=1.0\textwidth]{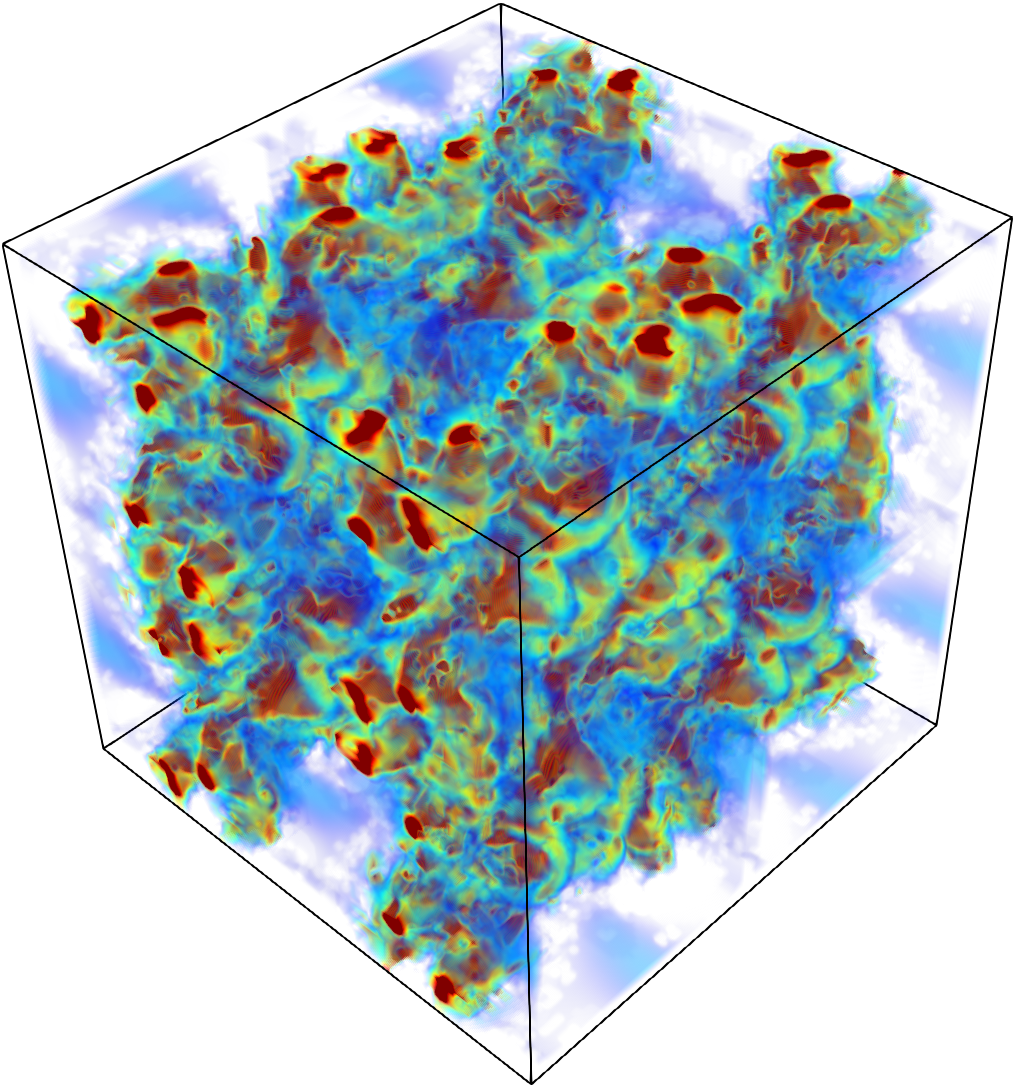}
        \caption{P4Q6 conformal, $t_c = 8$}
        \label{fig:conf64t8}
    \end{subfigure}
    \begin{subfigure}[b]{0.50\textwidth}
        \centering
        \includegraphics[width=1.0\textwidth]{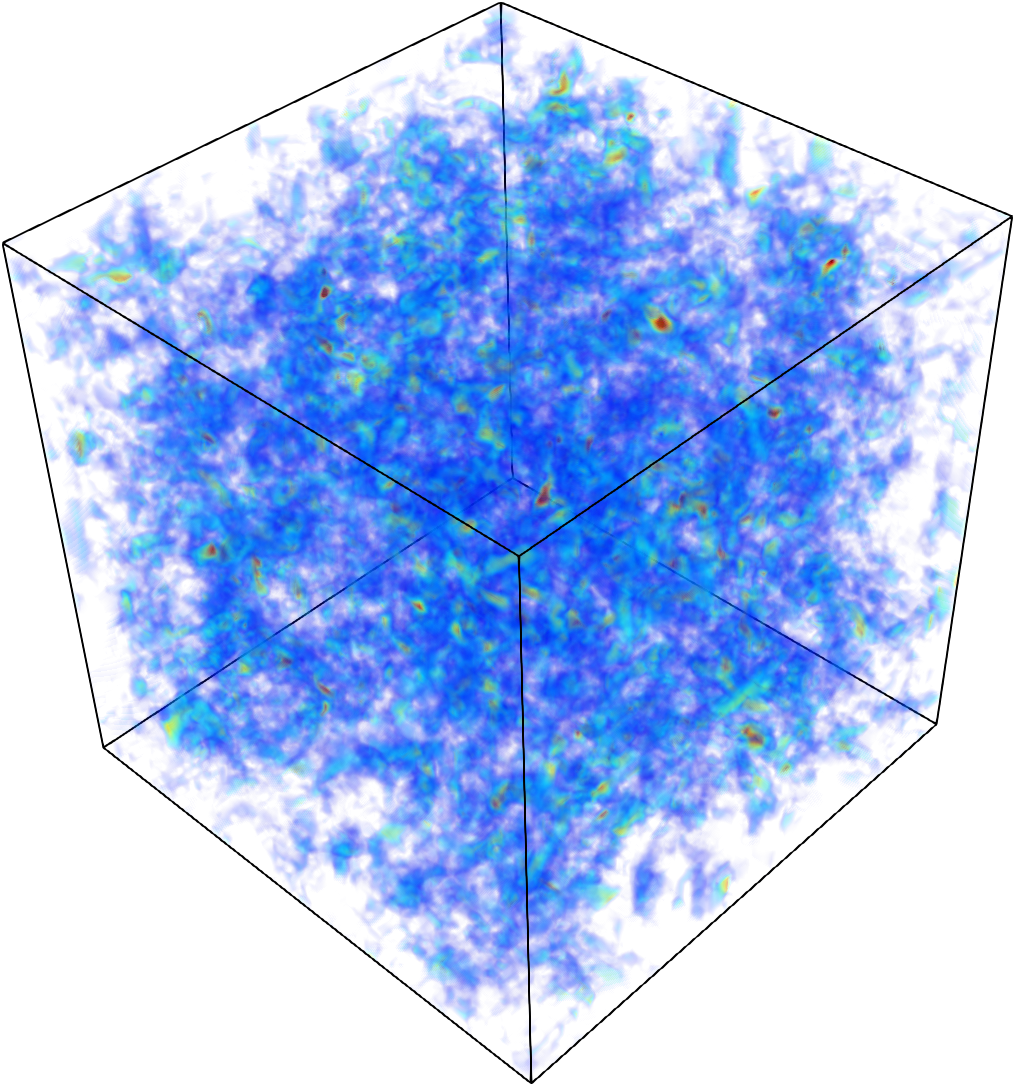}
        \caption{P4Q6 conformal, $t_c = 20$}
        \label{fig:conf64t20}
    \end{subfigure}
    \begin{subfigure}[b]{0.50\textwidth}
        \centering
        \includegraphics[width=1.0\textwidth]{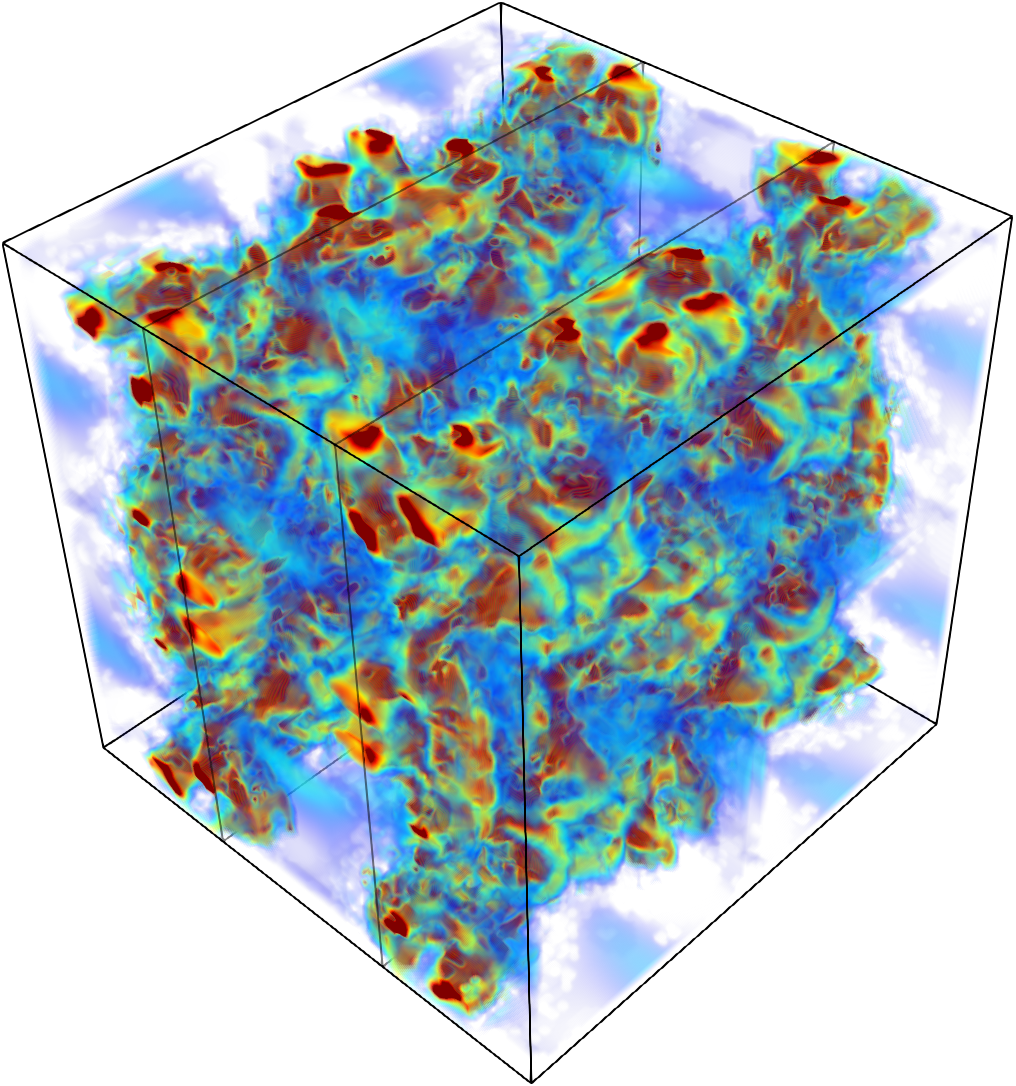}
        \caption{P4Q8 non-conformal, $t_c = 8$}
        \label{fig:nonconf64P4Q8t8}
    \end{subfigure}
        \begin{subfigure}[b]{0.50\textwidth}
        \centering
        \includegraphics[width=1.0\textwidth]{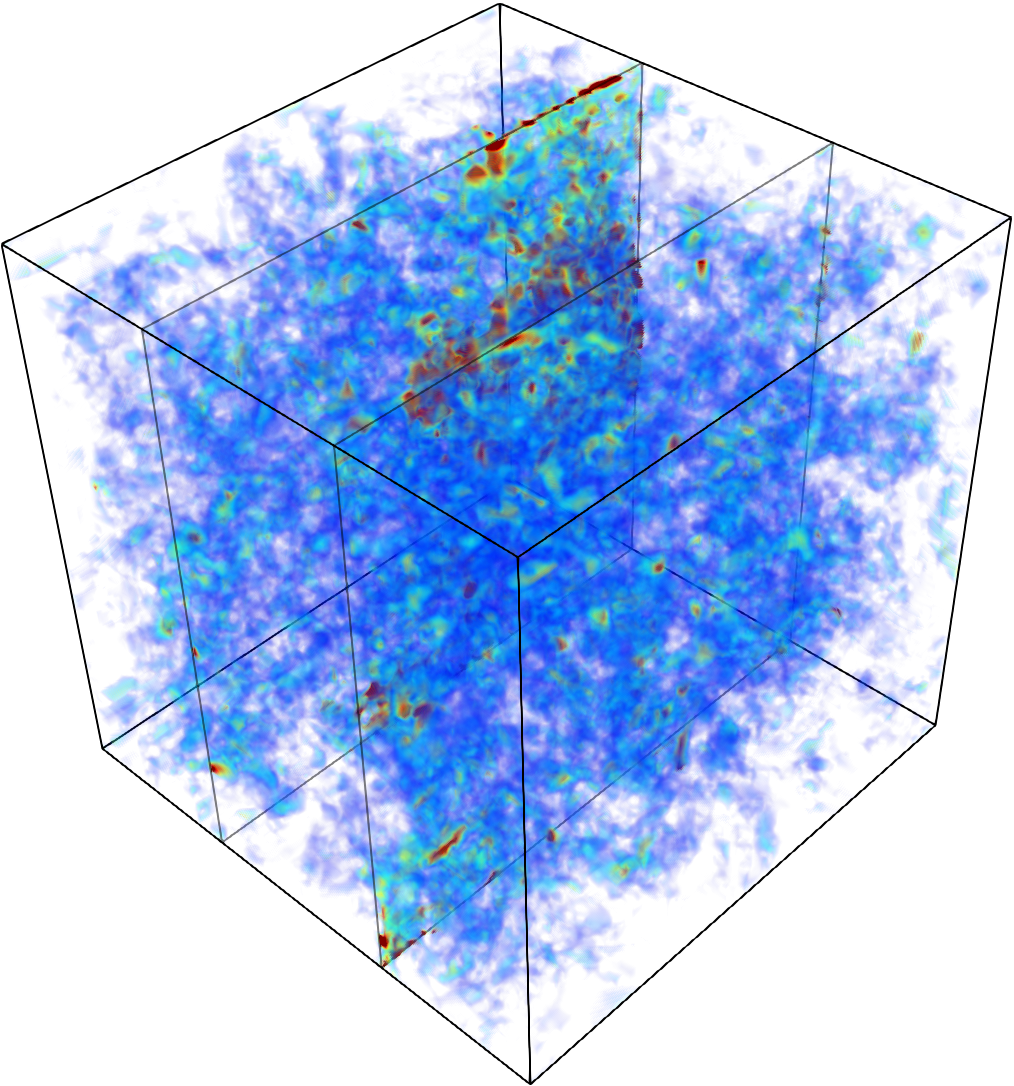}
        \caption{P4Q8 non-conformal, $t_c = 20$}
        \label{fig:nonconf64P4Q8t20}
    \end{subfigure}
    \begin{subfigure}[b]{0.50\textwidth}
        \centering
        \includegraphics[width=1.0\textwidth]{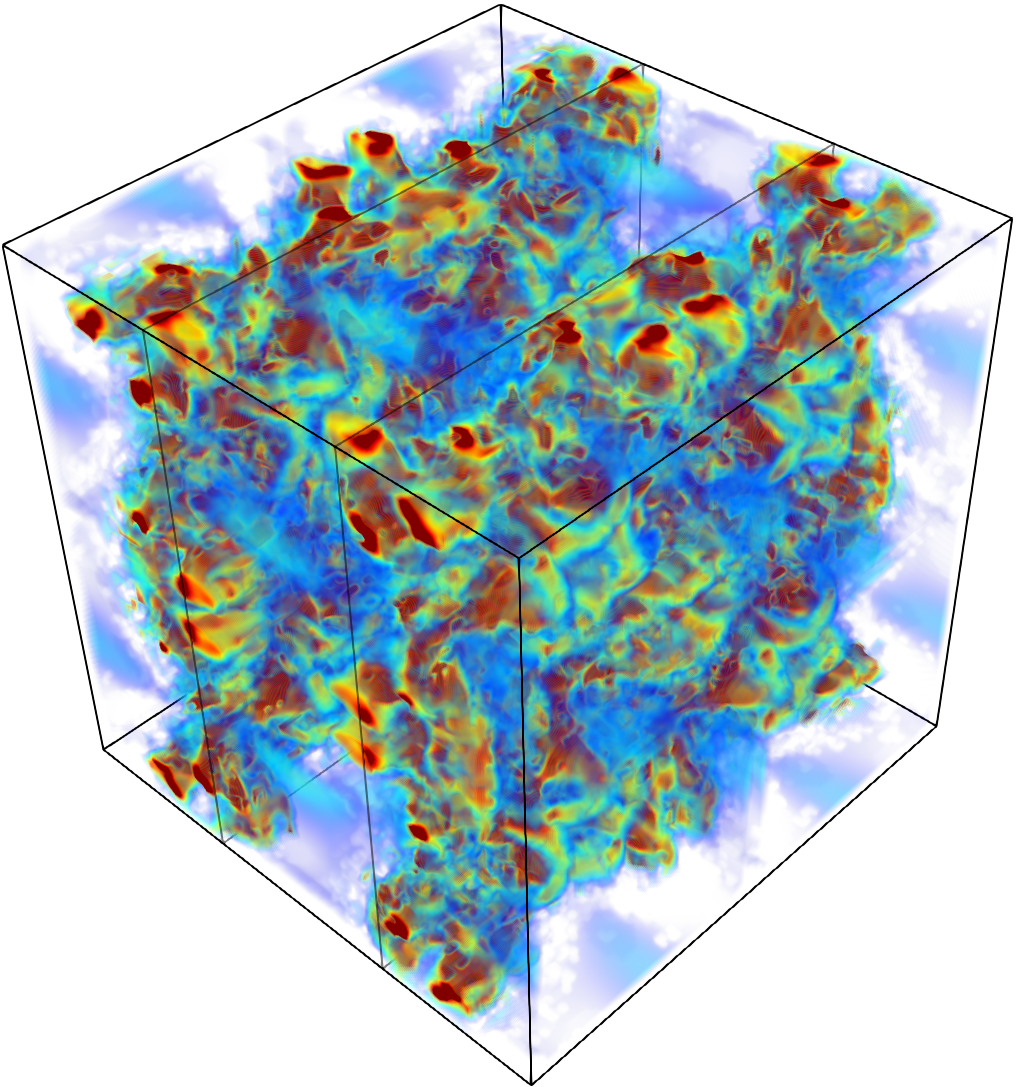}
        \caption{P4Q10 non-conformal, $t_c = 8$}
        \label{fig:nonconf64P4Q10t8}
    \end{subfigure}
        \begin{subfigure}[b]{0.50\textwidth}
        \centering
        \includegraphics[width=1.0\textwidth]{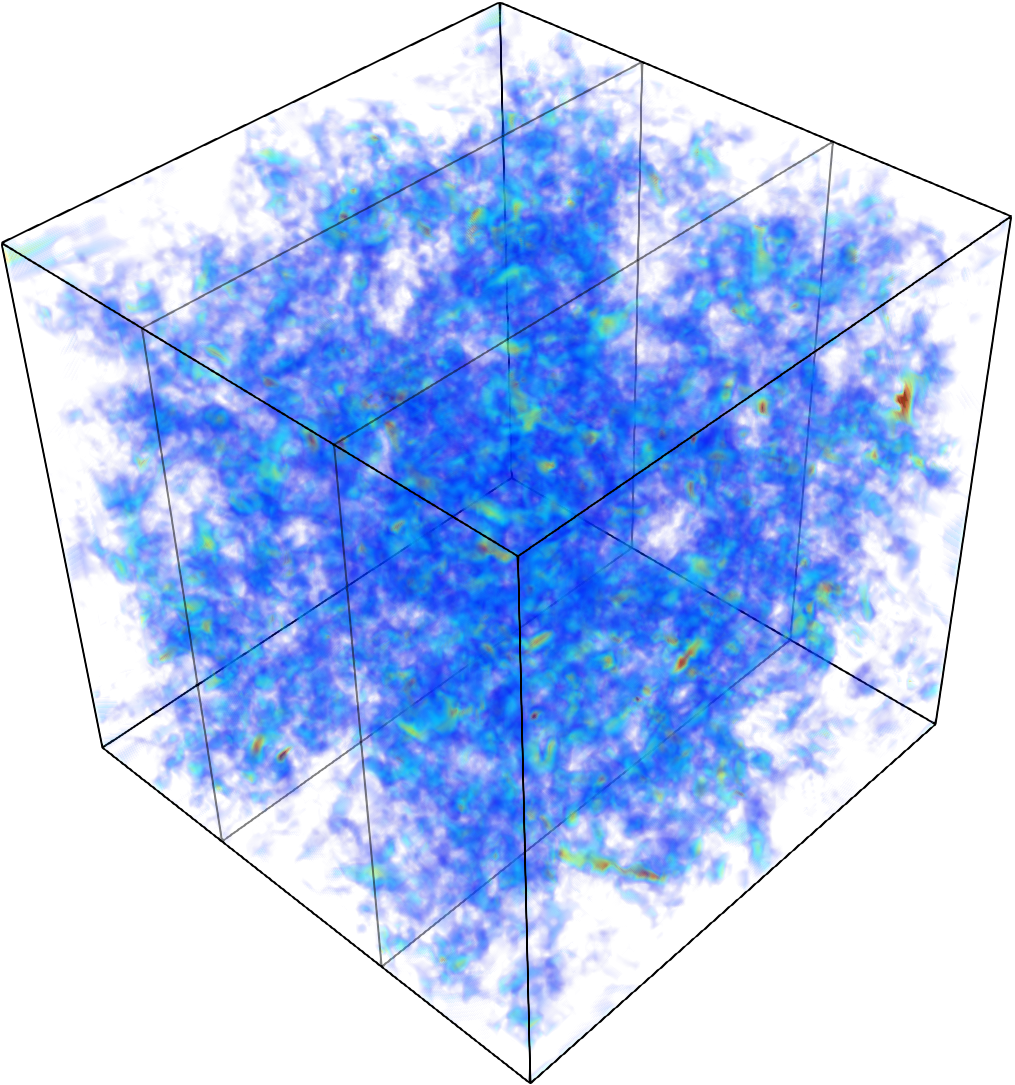}
        \caption{P4Q10 non-conformal, $t_c = 20$}
        \label{fig:nonconf64P4Q10t20}
    \end{subfigure}
    \caption{$64^3$ DOF cases showing $\|\bm{\omega}\|$ at $t_c = 8$ and $t_c = 20$.}
    \label{fig:vorticityNorm64}
\end{figure} 
\begin{figure}
        \begin{subfigure}[b]{0.50\textwidth}
        \centering
        \includegraphics[width=1.0\textwidth]{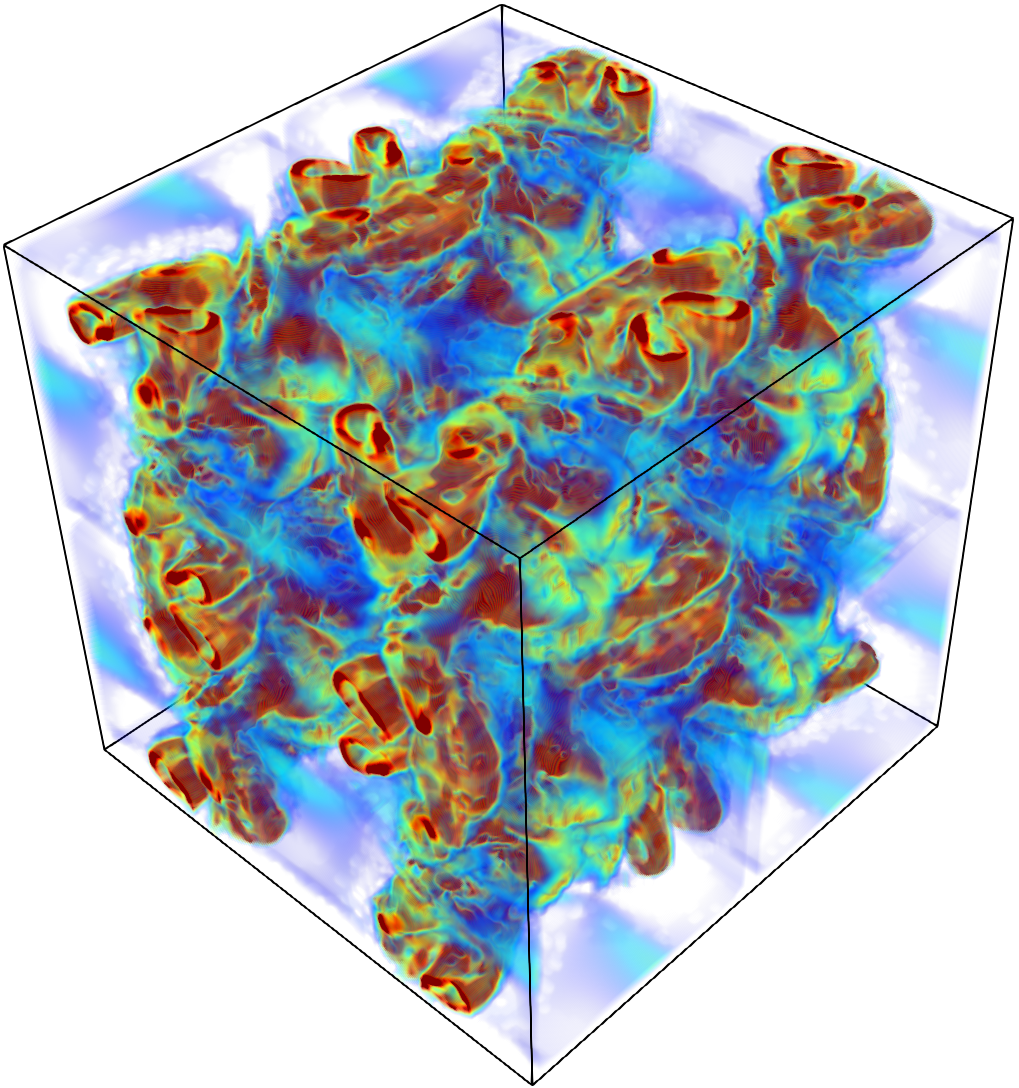}
        \caption{P4Q6 conformal, $t_c = 8$}
        \label{fig:conf128t8}
    \end{subfigure}
    \begin{subfigure}[b]{0.50\textwidth}
        \centering
        \includegraphics[width=1.0\textwidth]{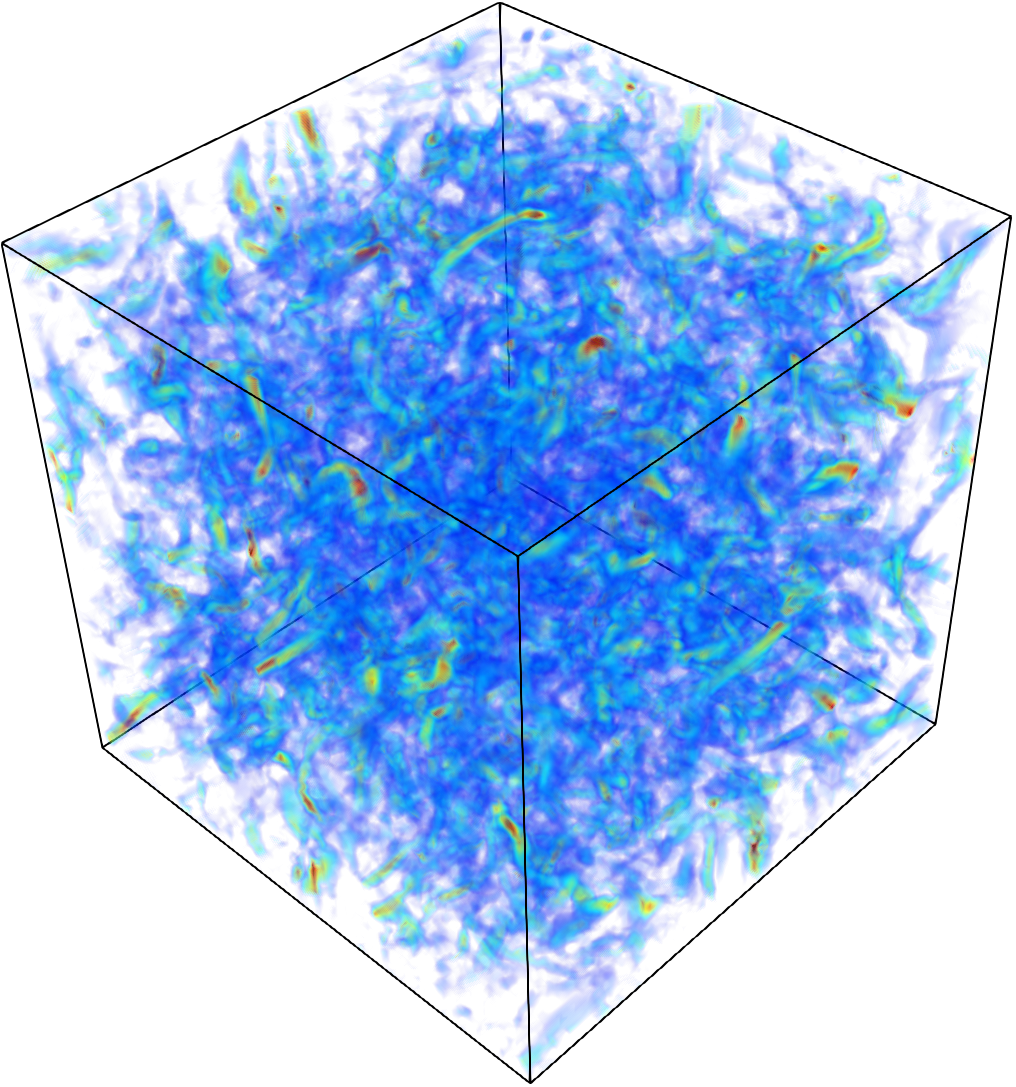}
        \caption{P4Q6 conformal, $t_c = 20$}
        \label{fig:conf128t20}
    \end{subfigure}
    \begin{subfigure}[b]{0.50\textwidth}
        \centering
        \includegraphics[width=1.0\textwidth]{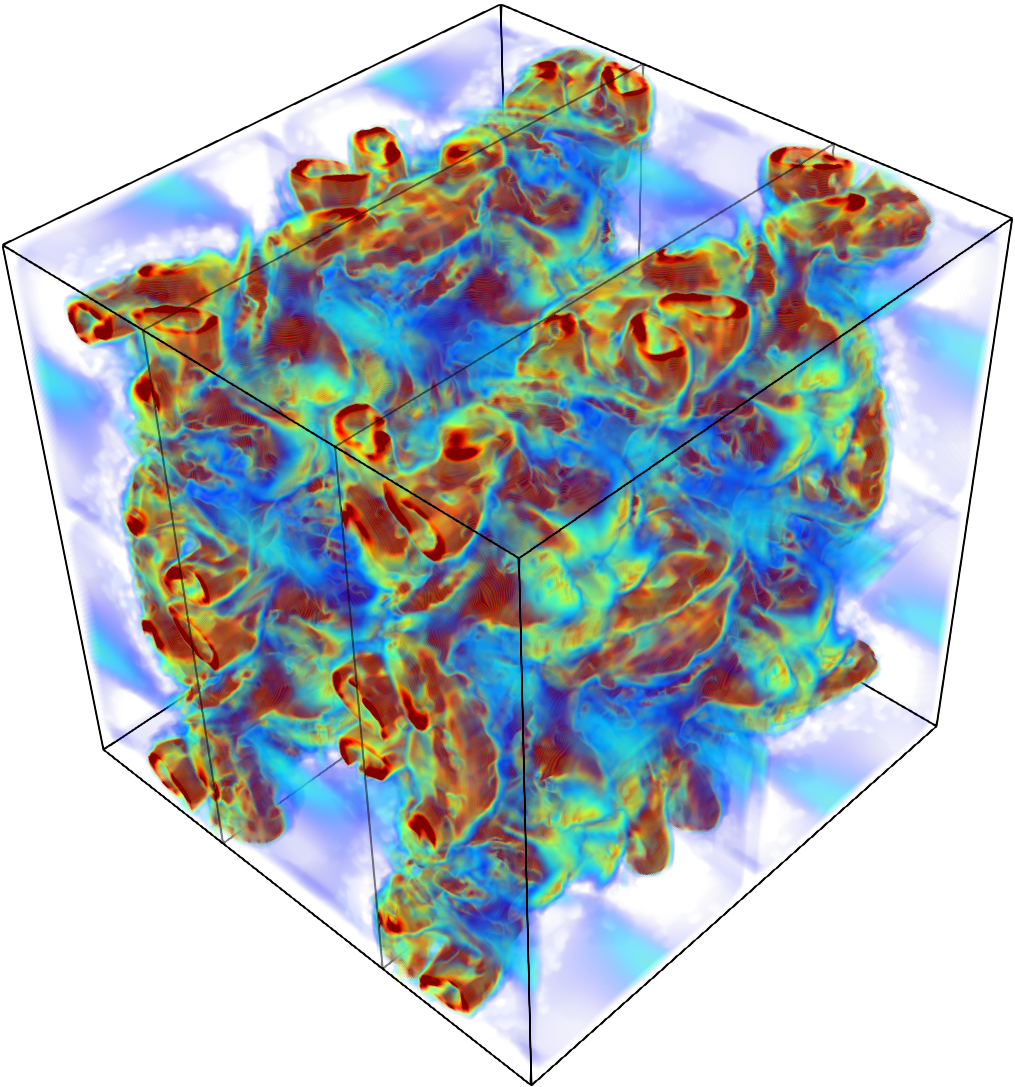}
        \caption{P4Q8 non-conformal, $t_c = 8$}
        \label{fig:nonconf128P4Q8t8}
    \end{subfigure}
        \begin{subfigure}[b]{0.50\textwidth}
        \centering
        \includegraphics[width=1.0\textwidth]{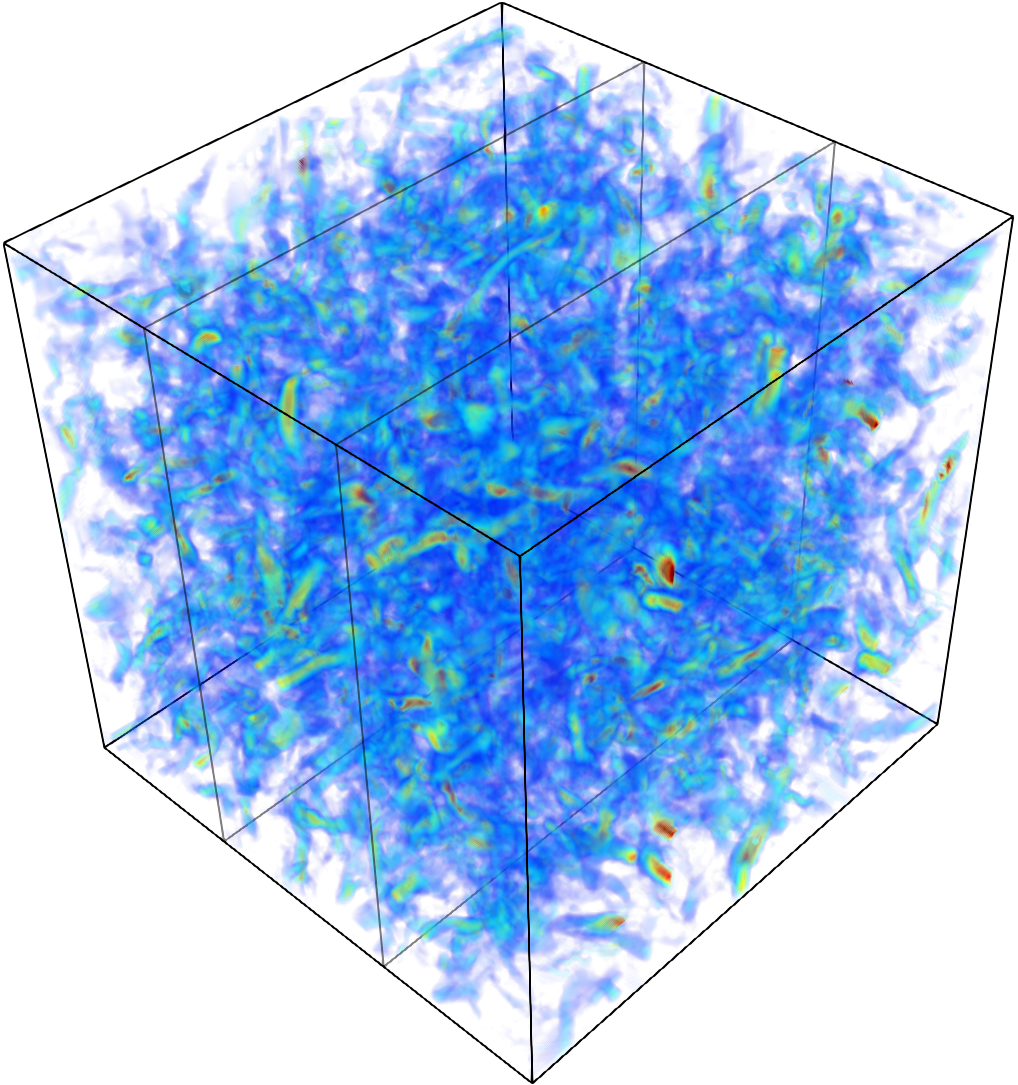}
        \caption{P4Q8 non-conformal, $t_c = 20$}
        \label{fig:nonconf128P4Q8t20}
    \end{subfigure}
    \caption{$128^3$ DOF cases showing $\|\bm{\omega}\|$ at $t_c = 8$ and $t_c = 20$.}
    \label{fig:vorticityNorm128}
\end{figure}

\begin{table}
\sisetup{round-mode=places, round-precision=4, table-format=1.4}
\centering
\caption{Computational costs for the Taylor-Green vortex cases.}
\label{tab:compCostsTGVortex}
\begin{tabular}{@{}lSSS@{}} \toprule
& \multicolumn{3}{c}{Avg. cost per timestep (s)} \\ \cmidrule(l){2-4} 
Case	      	    	        & {P4Q6} 		& {P4Q8}	    & {P4Q10}       \\ \midrule
Conformal -- $64^3$       		& 0.024675251   & 0.0384895	    & 0.06069	       	\\
Point-to-point -- $64^3$ 		& 0.0736983		& 0.108602091	& 0.151001274 	\\
Conformal -- $128^3$       	    & 0.191801619   & 0.335348	    & 0.69253	       	\\
Point-to-point -- $128^3$ 		& 0.303670069   & 0.529648609   & 1.09579  	    \\ \bottomrule 
\end{tabular}
\end{table}

\subsubsection{Examination of fluid structures}

Although the evolution of $\epsilon$ is an important global quantity of interest, we now consider snapshots of the vorticity norm $\|\bm{\omega}\|$ in order to view structures as they appear in the solution. Figure~\ref{fig:vorticityNorm64} for the $64^3$ DOF cases, and figure~\ref{fig:vorticityNorm128} for the $128^3$ DOF cases present volume renderings near the peak dissipation phase at $t_c = 8$ and at the final time $t_c=20$ when the smallest scale structures are present. The $64^3$ and $128^3$ DOF conformal cases are compared to the equivalent P4Q8 point-to-point cases, additionally a non-conformal $64^3$ P4Q10 is also shown, in order to examine the effects of higher dealiasing on the resulting solution. The interface locations have been depicted on the domain faces for reference in the non-conformal cases. It is apparent that at $t_c = 8$, very little difference can be seen between all cases and there is a close match between both conformal and non-conformal results. However, at $t_c = 20$ where the smallest vortex features appear, the $64^3$ DOF non-conformal P4Q8 case shown in figure \ref{fig:nonconf64P4Q8t20} shows a clear a buildup of vorticity that aligns with the non-conformal interface. By increasing the dealiasing to $Q=10$, this buildup is eliminated, as shown in figure~\ref{fig:nonconf64P4Q10t20}, and the resulting vorticity field closely resembles the conformal case in figure~\ref{fig:conf64t20}. At a higher resolution of $128^3$ DOF $\|\bm{\omega}\|$ snapshots for the same time points are shown in figure~\ref{fig:vorticityNorm128}. In this simulation it is clear that there is no accumulation of vorticity around the interfaces for the non-conformal P4Q8 case (fig.~\ref{fig:nonconf128P4Q8t20}) and the same small vortex features are present as can be seen in the conformal P4Q6 case (fig.~\ref{fig:conf128t20}).

\subsubsection{Computational cost}

The computational costs for the Taylor-Green vortex cases are shown in table~\ref{tab:compCostsTGVortex}. These simulations were run on 8$\times$ AMD Epyc 7742 ``Rome'' 64 core CPUs, hosted by the Isambard Tier 2 HPC facility, for a total of 512 cores. This shows the fairly large cost incurred by the dealiasing and the handling of the non-conformal interface, as well as the additional communication costs that are imposed in this setting. For example, going from the conformal $64^3$ P4Q6 case to the non-conformal $64^3$ P4Q10 case results in an approximately six times increase in average computational cost per timestep.

\section{Conclusions}
\label{sec:conclusions}

In this paper, we have compared the numerical performance of the point-to-point interpolation and mortar techniques, together with equivalent conformal cases, for a number of linear and non-linear hyperbolic conservation law problems. For problems that admit smooth solutions (i.e. which are adequately resolved in space), it is clear that either method is capable of performing equally well, both in terms of preserving the high-order convergence properties of the DG method, and also when considering the advection of structures across very long time periods. Likewise, when considering problems that are marginally- or under-resolved, it is equally clear that the mortar technique yields the most consistently accurate results when compared to the point-to-point interpolation approach.

Although there were relatively minor differences between the point-to-point and mortar methods for the linear Gaussian hump case in the presence of under-resolution, the isentropic vortex and Taylor-green vortex cases clearly highlight the care that must be taken when using the point-to-point method in such a regime. From the results we observe here, aliasing and oscillatory effects, owing to the discontinuity in polynomial interpolation across elemental interfaces, can have a significant impact on the ability of this method to accurately resolve flow features across long time periods or at a small scale. However, at the same time we note that it would be relatively unusual for higher-order fluid dynamics simulations to be performed in an implicit LES or under-resolved DNS regime without a significant level of dealiasing. As demonstrated in~\cite{Mengaldo2016}, running either compressible Euler or Navier-Stokes simulations without a comparable level of dealiasing to that we present here can yield inaccurate results and potentially lead to instability. Additionally, it is worth considering that in realistic fluid dynamics simulations of e.g. external aeronautics cases, most problems consist of inflow-outflow setups in which structures would be naturally removed from the domain within a far shorter time period than in the cases we consider here, which have been designed to deliberately test the numerical properties of each scheme.

We believe that there are three main factors to consider when choosing an interface handling technique for sliding or moving meshes: desired simulation accuracy, the capability to handle complex geometric interfaces and the resulting computational cost. In terms of accuracy, we have shown that the mortar method yields the results that are essentially identical to that of a conformal grid, and so for accuracy-critical simulations, this would certainly appear to be the most suitable strategy to adopt. Indeed in two dimensions, where the implementation is relatively straightforward, mortaring should be the first choice method to handle non-conformal grids. However the implementation challenge of constructing mortar elements across an arbitrary interface at high-order presents a significant obstacle in three dimensions. This makes the point-to-point method an attractive alternative, particularly in the context of highly parallel simulations. Results here show that the point-to-point interpolation method seems capable of handling a non-conformal interface in all cases, as long as appropriate precautions are taken by dealiasing to a sufficiently high level. The flexibility of supporting arbitrary interfaces in 3D, lower computational cost of interpolation on the interface and ease of implementation are advantages in this setting; however, they must be weighed against the major disadvantages of this method, which is the lack of a formal mass conservation and the requirement for dealiasing (which then further increases computational cost).

Further investigation is warranted to investigate the minimum amount of dealiasing required for the point-to-point interpolation method, and whether it is suitable to overintegrate only on the interface skeleton elements, which would undoubtedly significantly reduce the overall computational cost. Another aspect that we do not consider in terms of computational cost is the relative effort required to set up a moving grid, where the creation of mortars and evaluation of interpolation points needs to be performed at every timestep. It would also be beneficial to compare an efficient 3D mortar method implementation and the point-to-point interpolation method with conformal solutions in a 3D setting for more complex flow problems. Additionally, investigating problems involving shocks would be another area of research to pursue, to validate the point-to-point interpolation and mortar method under more demanding transonic or supersonic conditions.

\subsection*{Acknowledgements}		
DM acknowledges support from the EPSRC Platform Grant PRISM under grant EP/R029423/1 and the ELEMENT project under grant EP/V001345/1. This work used the Isambard UK National Tier-2 HPC Service (\url{http://gw4.ac.uk/isambard/}) operated by GW4 and the UK Met Office, and funded by EPSRC under grant EP/P020224/1.
\newpage
\bibliography{elsarticle-template}

\end{document}